\documentclass[12pt]{article}
\usepackage{graphicx}
\usepackage{psfrag}
\usepackage{amsmath}
\usepackage{amssymb}
\usepackage{amsthm}

\newtheorem{definition}{Definition}
\newtheorem{corollary}{Corollary}
\newtheorem{theorem}{Theorem}
\newtheorem{remark}{Remark}
\newtheorem{hypothesis}{Hypothesis}
\newtheorem{proposition}{Proposition}
\newtheorem{example}{Example}

\begin{document}

\title{Period-doubling cascades galore}

\author{Evelyn Sander and James A. Yorke}

\maketitle

\begin{abstract}
  The appearance of numerous period-doubling cascades is among the
  most prominent features of {\bf parametrized maps}, that is, smooth
  one-parameter families of maps $F:R \times {\mathfrak M} \to
  {\mathfrak M}$, where ${\mathfrak M}$ is a smooth locally compact
  manifold without boundary, typically $R^N$. Each cascade has
  infinitely many period-doubling bifurcations, and it is typical to
  observe -- such as in all the examples we investigate here -- that
  whenever there are any cascades, there are infinitely many cascades. We
  develop a general theory of cascades for generic $F$. We illustrate
  this theory with several examples. We show that there is a close
  connection between the transition through infinitely many cascades
  and the creation of a horseshoe.

\end{abstract}

\section{Introduction}
\begin{figure}
\begin{center}
\includegraphics[width=\textwidth]{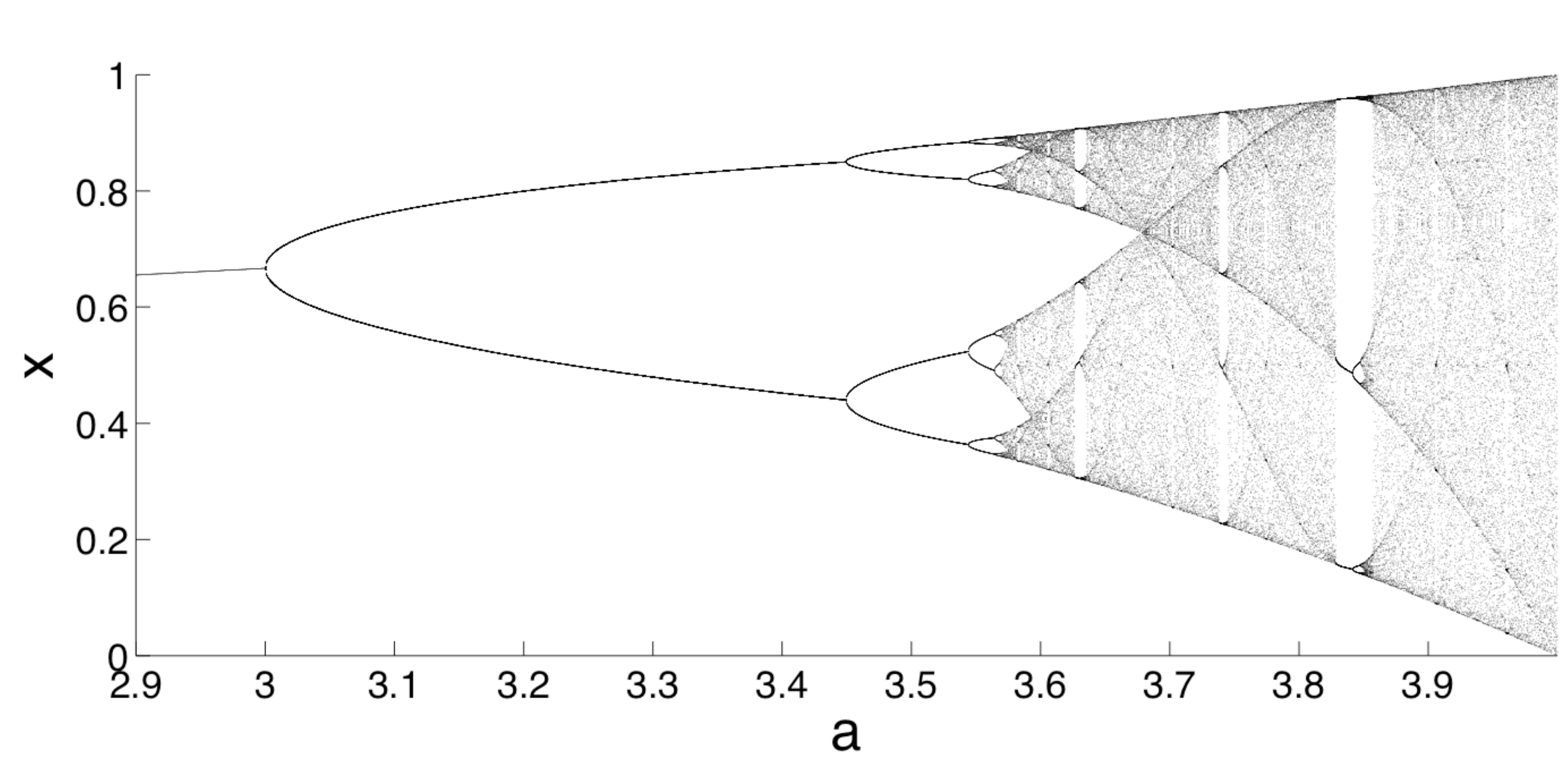}
\caption{\label{fig:typicalquadratic} The attracting set for the
  logistic map: $F(a,x)=a x (1-x)$. That is, for each fixed parameter value, the attracting set in $[0,1]$ is shown. There are infinitely
  many cascades of attractors. This is the bifurcation diagram most
  frequently displayed to illustrate the phenomenon of period-doubling
  cascades. However, cascades occur for much more complex dynamical systems that are completely unrelated to quadratic maps.}
\end{center}
\end{figure}

\begin{figure}
\begin{center}
\includegraphics[width=\textwidth]{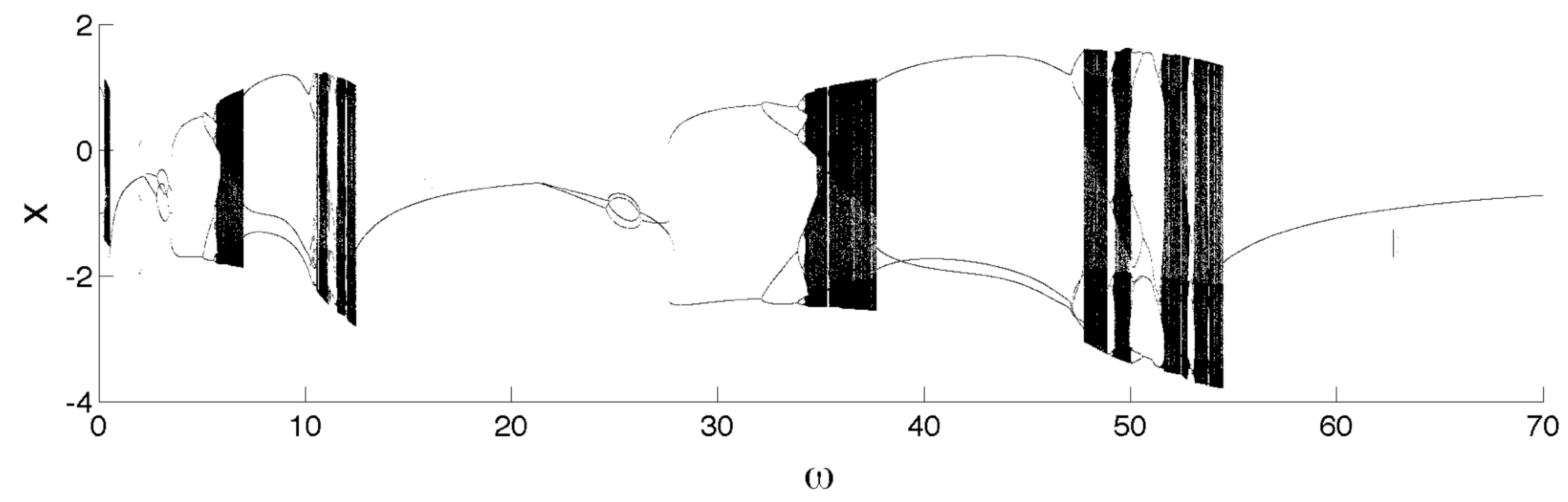}
\caption{\label{fig:duffing} The attracting set for the double-well
  Duffing equation: 
   $ u''(t) + 0.3 u'(t) - u(t) + (u(t))^3 +0.01 =  \omega \sin t.$
This equation is periodically forced with period
  $2\pi$. Therefore the time-$2\pi$ map is a diffeomorphism on $R^2$ parametrized by $\omega$.  Depicted here is the
  attracting set of $F$, projected to the $ ( \omega,u'(t) ) $-plane.  
The constant $0.01$ has been added to destroy symmetry in order to avoid nongeneric symmetry-breaking bifurcations. }
\end{center}
\end{figure}

A major goal in dynamical systems is to better explain what is seen.
There are many reports of numerically observed (period-doubling) cascades in the
contexts of  maps, ordinary differential equations,
partial differential equations, and even in physical experiments.  In
each of these contexts, whenever one cascade is seen, an infinite
number are observed. There is also considerable numerical evidence
that cascades occur as a dynamical system becomes more chaotic. The
concept of cascades is usually associated with the orderly picture of
the attractor structure of a quadratic map depicted in
Figure~\ref{fig:typicalquadratic}. The structure of quadratic maps is
well understood, and it is not hard in this context to rigorously show
that cascades exist. However, most dynamical systems have bifurcation
structures which are a good deal more complicated than that of a
quadratic map (cf. the double-well Duffing equation in
Figure~\ref{fig:duffing}), implying that the simple explanation for
cascades of the quadratic map does not generalize. In this paper, we
develop a series of general criteria for the existence of cascades in
the context of parametrized maps.  We explain why cascades occur with
infinite multiplicity. Our results give a rigorous explanation for the
link between cascades and chaos.  The method of approach lays a
general framework for understanding cascades, even for observable
systems for which the underlying model equation is unknown.

We now give a heuristic description of our results.

We call a period-$k$ orbit of a parametrized
map $F:R \times {\mathfrak M} \to {\mathfrak M}$ a {\bf flip} orbit if
its Jacobian matrix $D_xF^k(\lambda, x)$ has an odd number of eigenvalues
less than -1, and -1 is not an eigenvalue.  Otherwise the orbit is
{\bf nonflip}.  We denote the space of nonflip orbits of $F$ in
$R\times {\mathfrak M}$ under the Hausdorff metric by  
$PO_{nonflip}(F)$.  We restrict our attention to a specific residual
set  (which will be precisely defined) of the $C^\infty$ maps $F:R\times {\mathfrak M} \to {\mathfrak
  M}$. We call maps in this residual set  {\bf generic}. All our results are
for generic maps. 

{\bf Components of $PO_{nonflip}(F)$ are one-manifolds for generic $F$}. 
For a specific residual set of $F$, which we call {\em generic $F$}, we show in Theorem~\ref{thm:OneManifold} ($M_1$) that  all
connected components of $PO_{nonflip}(F)$ are one-manifolds; that is,
they are either homeomorphic to circles or to open
intervals. Furthermore, the ratio of the periods of two orbits in a
component of $PO_{nonflip}(F)$ is always a power of 2.
This result may seem counterintuitive to those familiar with
cases where periodic points are dense in a space. But our result is
about orbits, not points. This result is one of the keys to
understanding cascades. From now on, we use the term {\bf component}
to denote a connected component of $PO_{nonflip}(F)$. When a component
is homeomorphic to an open interval, we call it an {\bf open arc}.

{\bf Definition of Cascade}.  Cascades were first reported by Myrberg
in 1962~\cite{myrberg:62}.  See also Robert May~\cite{may:74}.  We
define a {\bf cascade} as the following type of
subarc of an open arc $A$: Let $k$ denote the smallest period of the
orbits in $A$.  A cascade is a half-open subarc that contains orbits with all of
the periods $k, 2k, 4k, 8k, \cdots$ such that it contains precisely
one orbit of period $k$.  We refer to such a cascade as a {\bf period-$k$
  cascade}. A cascade is homeomorphic to $[0,1)$ where $0$ maps to
the single period-$k$ orbit.

\begin{figure}
\begin{center}
\includegraphics[width=.45\textwidth]{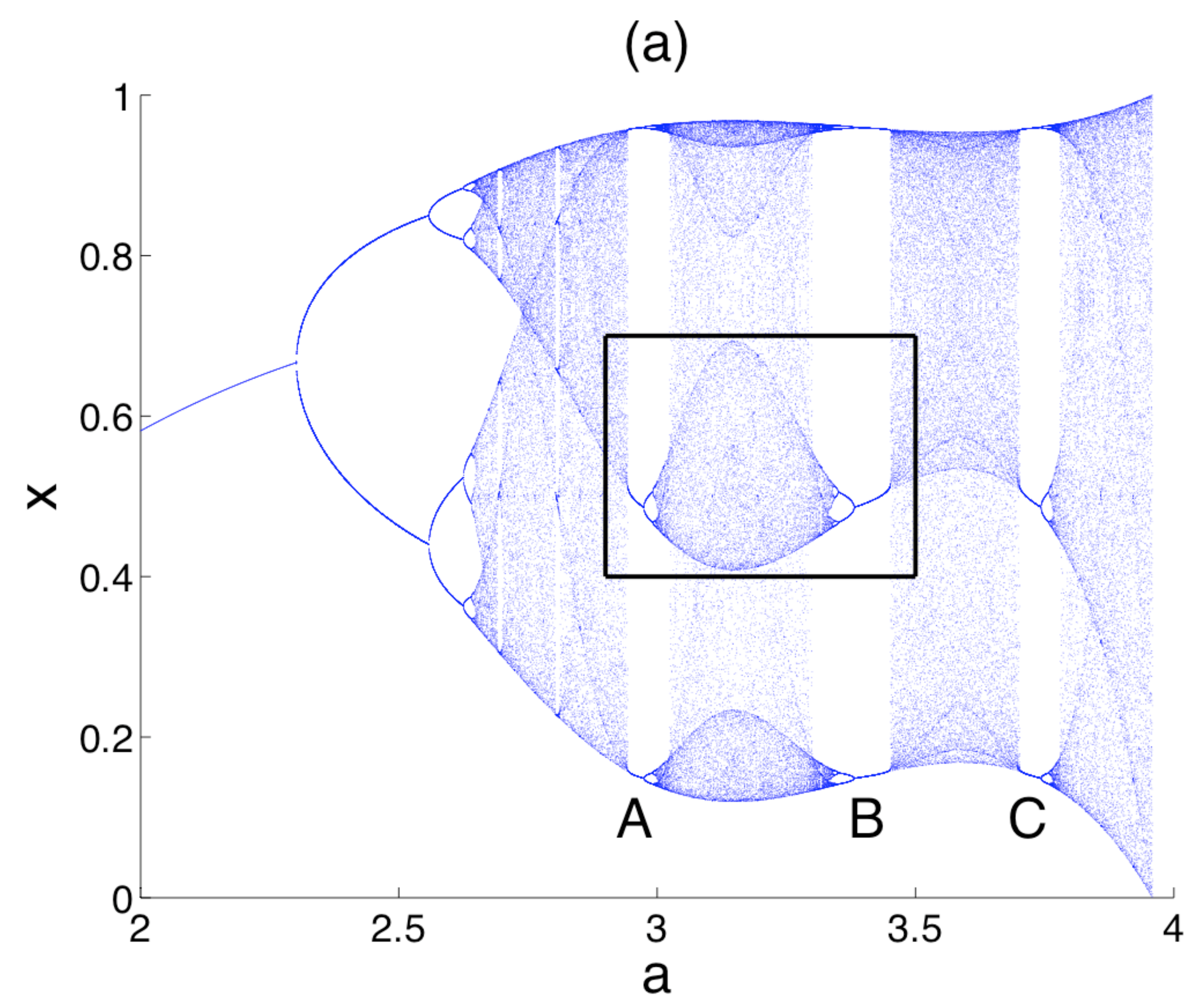}
\includegraphics[width=.45\textwidth]{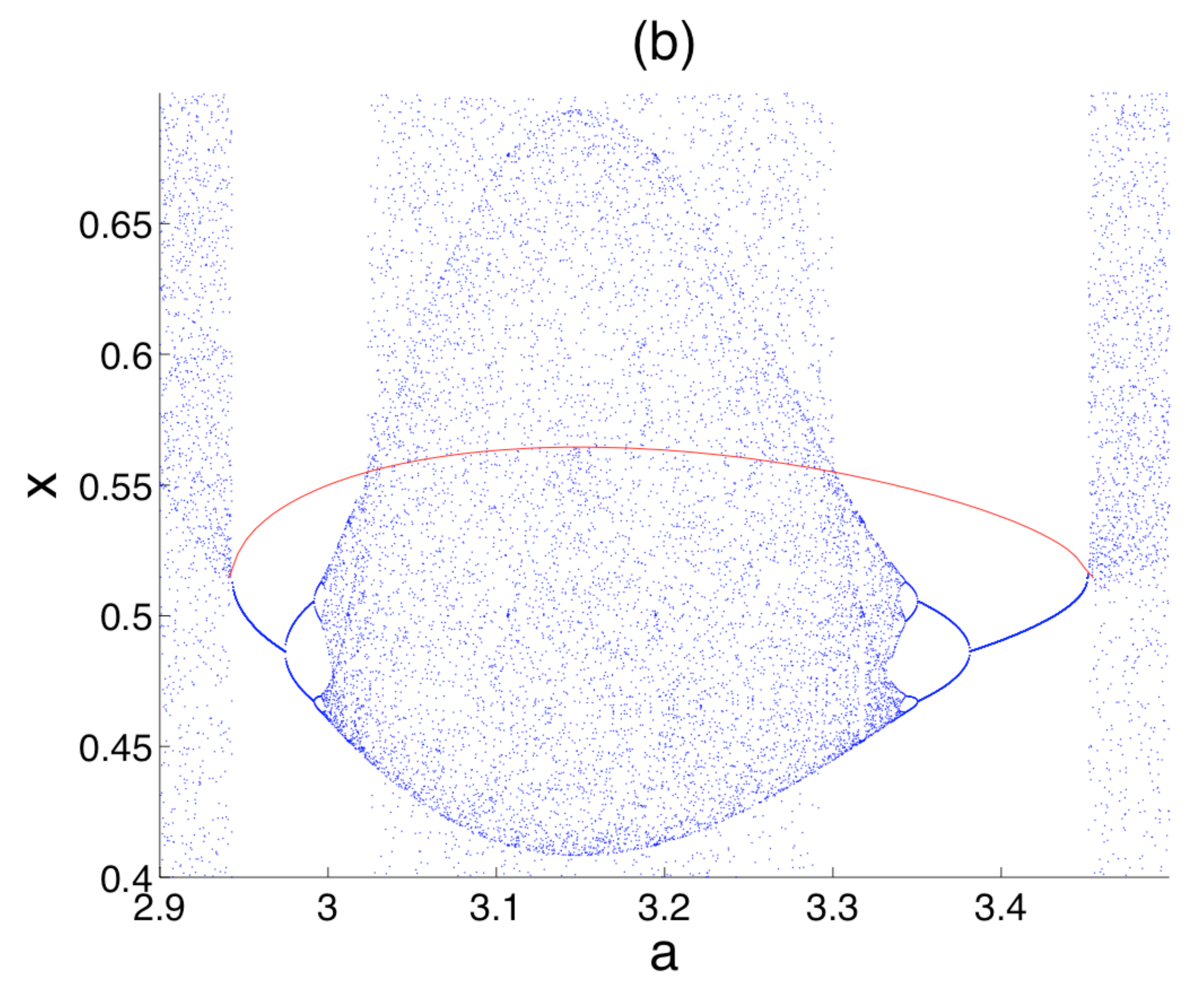}
\caption{\label{fig:boundedcascades} Bounded cascades come in pairs. This figure
  shows the numerically computed attracting set for a modified
  logistic map $F(a,x)= h(a) x (1-x) $, where $h(a)=
  a (1.18 + 0.17 \cos(2.4 a))$. Here $h$ is not monotonic. As the parameter $a$ increases, $h(a)$ increases, then decreases, and again increases.
Hence $h(a)$ passes three times over the region where the logistic map has a period-three cascade. Figure (a) displays three
  ``windows'' (the parameter ranges labeled A,B, and C) of
  period-doubling sequences starting from period-three orbits. 
When parameter $a$ yields a one-piece chaotic attractor,
  the upper edge of the chaotic set is the image of $x=1/2$,
  corresponding to h(a)/4, the maximum value of $F(a,\cdot)$. Figure (b) shows a blowup of the box
  in Figure (a). We have added (in red) a curve of unstable nonflip period-three points
revealing that the two cascades are in the same component.  Hence we see three period-three cascades, revealing a phenomenon that often happens in much more complicated systems. 
Namely,  that bounded
  cascades can be created or destroyed in {\it pairs that are in the same component}. 
 }
\end{center}
\end{figure}

We call an open arc $A$ {\bf bounded} if there is a compact subset
of $R \times {\mathfrak M}$ that contains all the orbits of $A$.
Otherwise $A$ is {\bf unbounded}. If a cascade is contained in an
unbounded open arc, we call it {\bf unbounded}; otherwise we call it
{\bf bounded}. We show that if a component is a bounded open arc, then
it always contains two disjoint cascades, and we refer to these as
{\bf paired cascades}. See Figure~\ref{fig:boundedcascades}.

{\bf The orbit index.} 
An essential tool for studying cascades is the {\bf orbit index}, a
topological index of a periodic orbit taking on a value in $\{ 0, -1,
+1 \}$. A periodic orbit in $PO_{nonflip}(F)$, has an orbit index of
either $-1$ or $+1$.  In Theorem~\ref{thm:OneManifold} $(M_2)$,
we state and prove that each component has a preferred orientation that can be
determined at each hyperbolic orbit by computing its orbit index.

Section~\ref{sec:abstracttheorem} contains Theorem~\ref{thm:main},
establishing existence of cascades in a bounded parameter region based on a
small amount of information about the types of periodic orbits on the
boundary. The result is established by using the preferred orientation via the orbit index
to show that an open arc in a bounded parameter region is bounded and therefore contains a cascade.
In order to illustrate the implications of these abstract statements, we sketch some
prototypical examples. See the relevant applications (labeled
corollaries) for the precise statements.

\bigskip 

{\bf Corollary~\ref{cor:SmaleHorseshoe}: Maps with Smale
horseshoes in dimension 2.}   
Assume ${\mathfrak M}$ is a two-manifold and $F$ is generic. Assume there are parameter
values $\lambda_0$ and $\lambda_1$ for which $F(\lambda_0, \cdot)$ has
at most finitely many nonflip orbits, and $F(\lambda_1, \cdot)$ has infinitely many
nonflip saddles and at most a finite number of nonflip attractors and repellers. Then $F$ has infinitely many cascades between $\lambda_0$ and $\lambda_1$.

This result follows from the fact that all but a finite number of the
nonflip orbits in the horseshoe are in components that contain either
one or two cascades, and all but a finite number of the distinct
nonflip orbits in the horseshoe are in different components.

\bigskip 

{\bf Corollary~\ref{cor:OffOnOffChaos}: Maps with ``Geometric 
Off-On-Off Chaos'' in dimension 2.} 
These are maps that first have no chaos, then chaos appears, and then it disappears as $\lambda$ is varied. This scenario appears to happen with the time-$2\pi$ maps of the single-well and double-well Duffing and the forced damped pendulum. 

Assume that Corollary~\ref{cor:SmaleHorseshoe} above applies for the interval $[\Lambda_1,\Lambda_2]$.   
Assume there is parameter $\Lambda_3$ with
$\Lambda_1 < \Lambda_2 < \Lambda_3$ for which 
$F(\Lambda_3,\cdot)$ also has at most a finite number of nonflip orbits, 
Then $F$ has infinitely many bounded paired cascades, and it has at most finitely many unbounded
cascades that have any orbits with $\lambda$ values in $[\Lambda_1,
\Lambda_3]$.

\bigskip
{\bf Corollary~\ref{cor:1D-quadratic}: Large-scale perturbation of quadratic maps.}
Consider the parametrized one-dimensional map
\[ F(\lambda, x) = \lambda - x^2 + g(\lambda, x), \] where $\lambda$
and $x \in R$, and $g: R \times R \to R$.  For each generic, smooth,
$C^1$-bounded function $g$, the map $F$ has exactly the same number of
unbounded period-$k$ cascades as occur in the case $g\equiv 0$, and we
give a recursive formula for that number. $F$ may have extra cascades,
but they are all bounded paired cascades. That is, the number of
unbounded cascades is {\em robust under large-scale perturbations.}

\bigskip

\begin{figure}[t]
\begin{center}
\includegraphics[width=.45\textwidth]{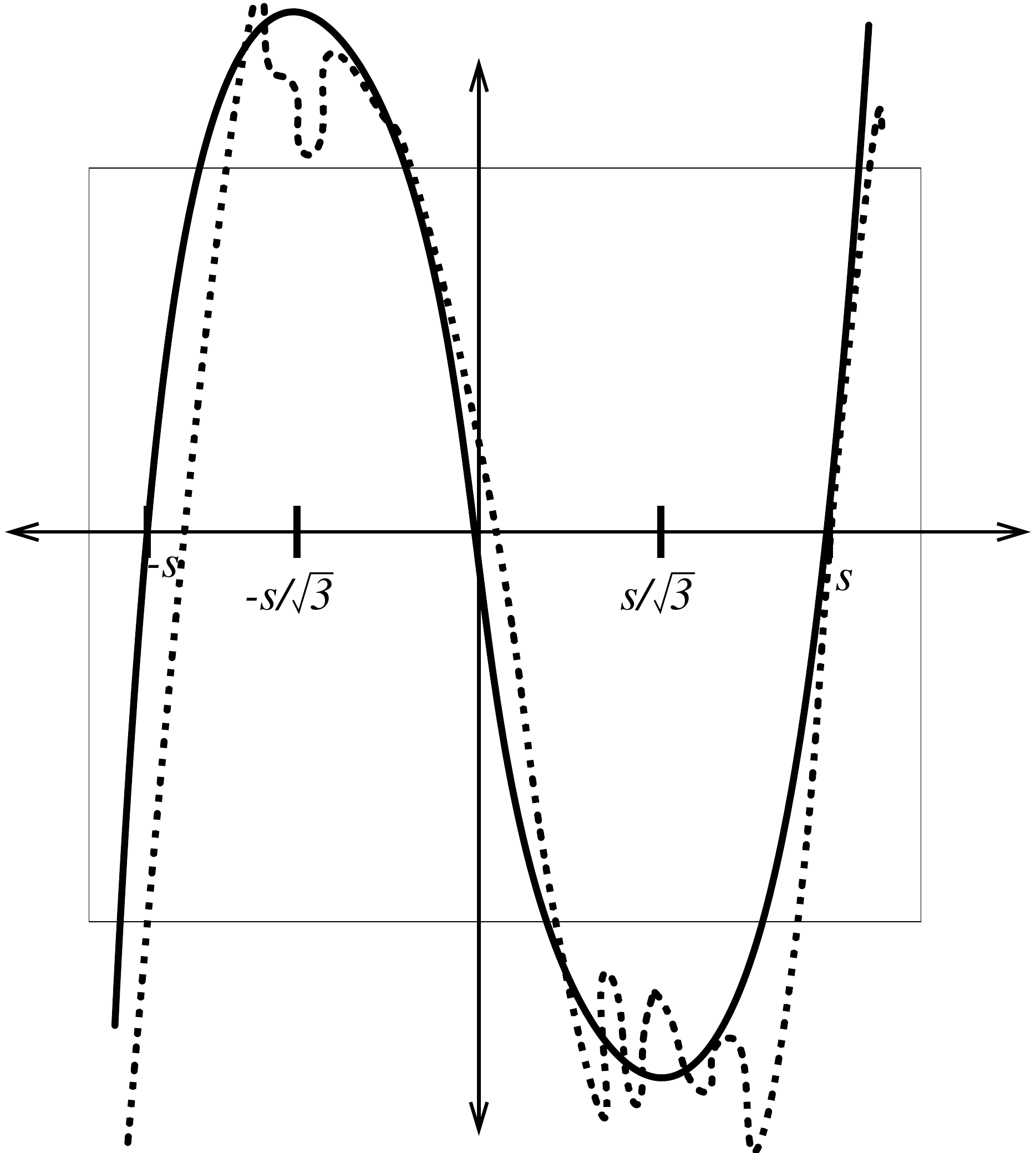}
\caption{\label{fig:cubicheur} 
A large-scale perturbation of the cubic parametrized map for large $\lambda$: This figure is a heuristic depiction (not drawn to scale) 
 of (solid line) the standard cubic, $x^3-\lambda_1 x$, for large
 $\lambda_1$, along with (dotted line) an example of a large-scale perturbation
 of the
 cubic $F(\lambda_1,x)$ as in Corollary~\ref{cor:cubic}. In the proof
 of the corollary, we establish and use facts depicted here. Namely that,
 letting $s=\sqrt{\lambda_1}$, the standard cubic has zeroes at $0$
 and $\pm s$, and has critical points at $\pm
 s/\sqrt{3}$. Furthermore, it has derivative with norm greater than
 one inside the box $[-2s,2s]^2$, and the graph of both the critical
 points and of the points near the edge of the box are outside of the
 box. } 
\end{center}
\end{figure}

{\bf Corollary~\ref{cor:cubic}: Large-scale perturbations of the cubic map.} The results in Corollary~\ref{cor:1D-quadratic} are not related to whether the base map is quadratic.
Similar behavior occurs for the parametrized cubic map 
\[ x^3-\lambda x +g(\lambda,x). \]
This parametrized map has infinitely many cascades for a residual set of $C^\infty$ functions $g$ that are $C^1$ bounded. The number of period-$k$ cascades differs from the quadratic case. This is to be expected since the behavior of the cubic map for large $\lambda$ reflects a shift on three symbols and so is more complex than that of the tent map. See Figure~\ref{fig:cubicheur}.

\begin{figure}
\begin{center}
\includegraphics[width=.45\textwidth]{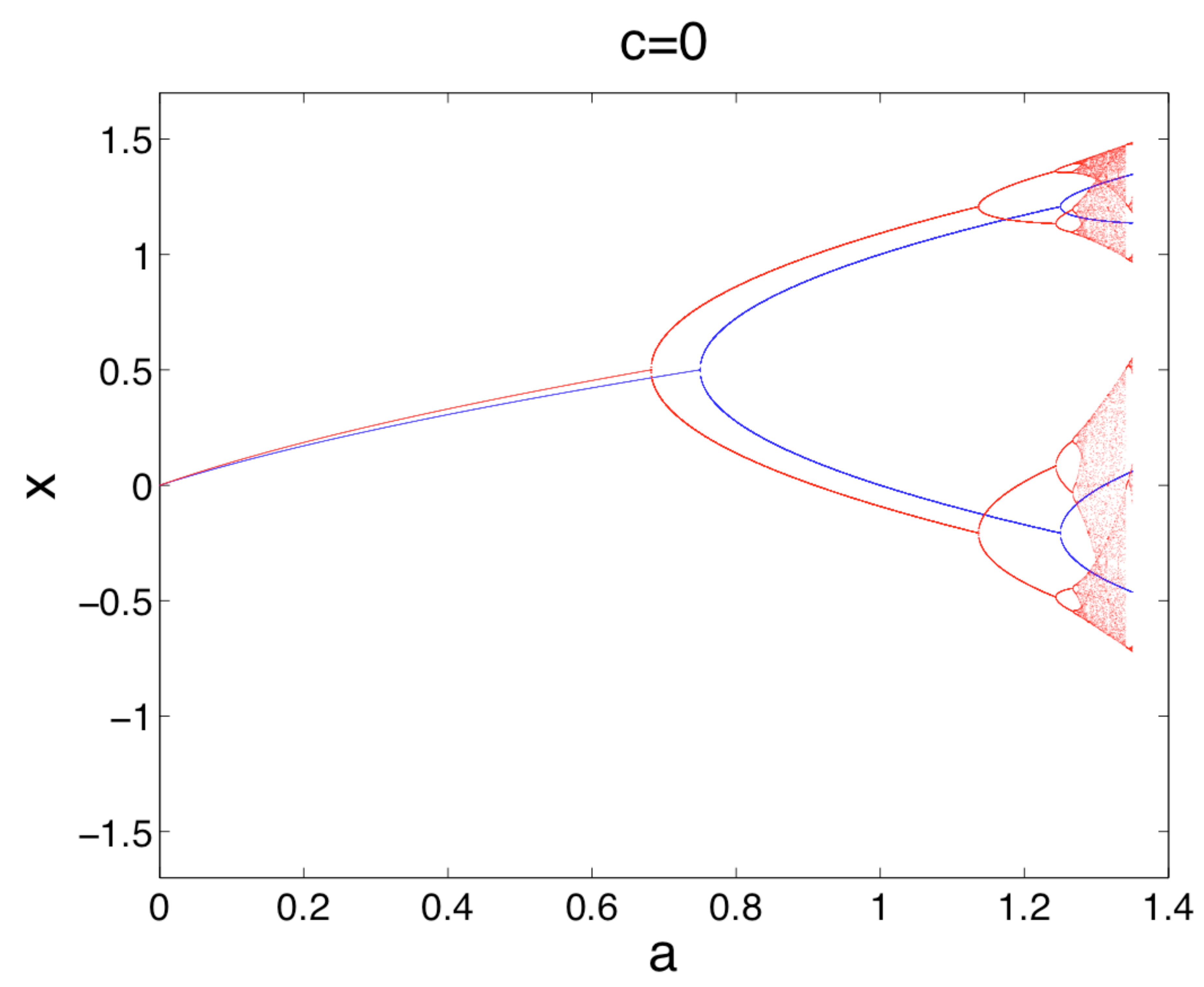}
\includegraphics[width=.45\textwidth]{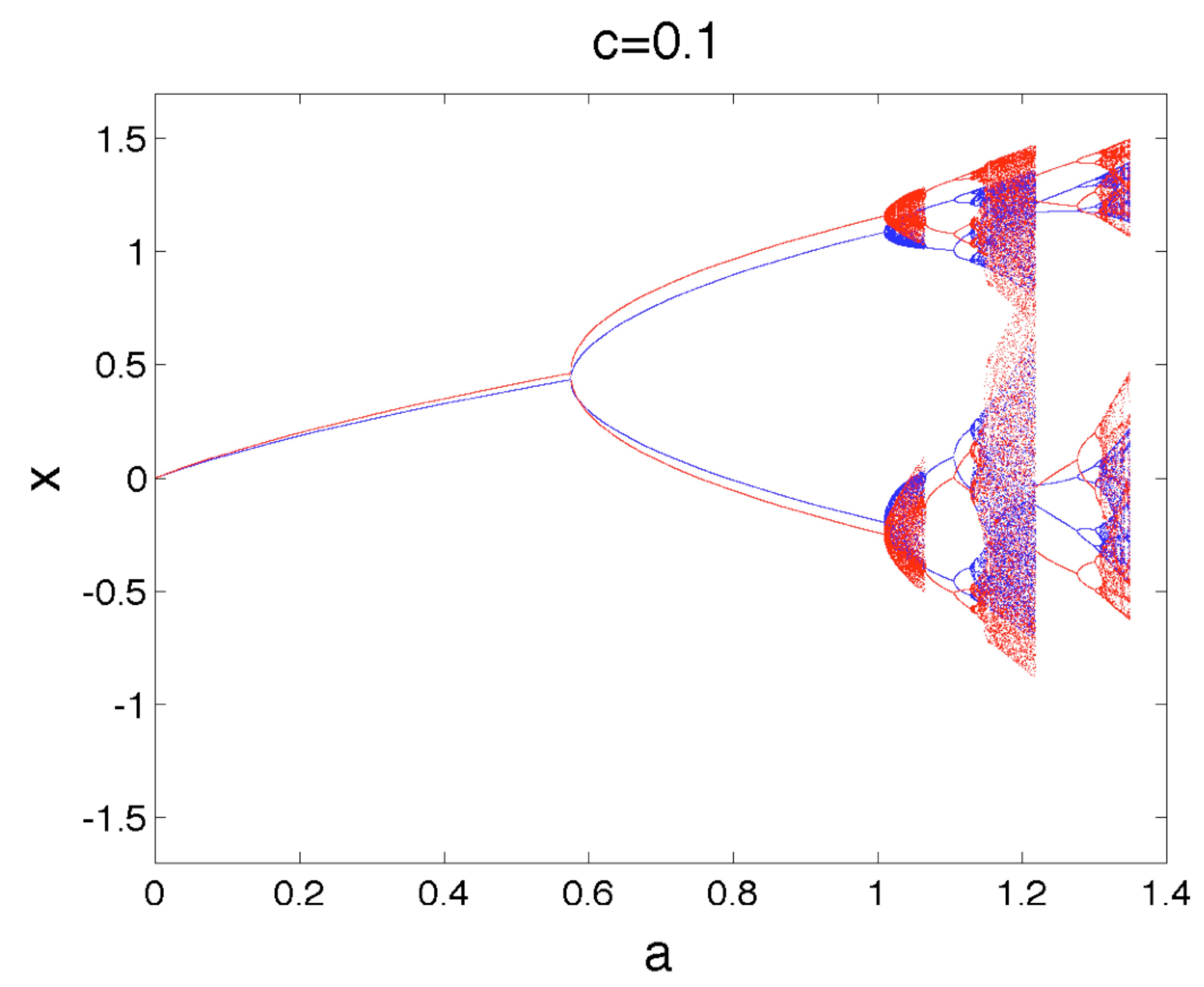}\\
\includegraphics[width=.45\textwidth]{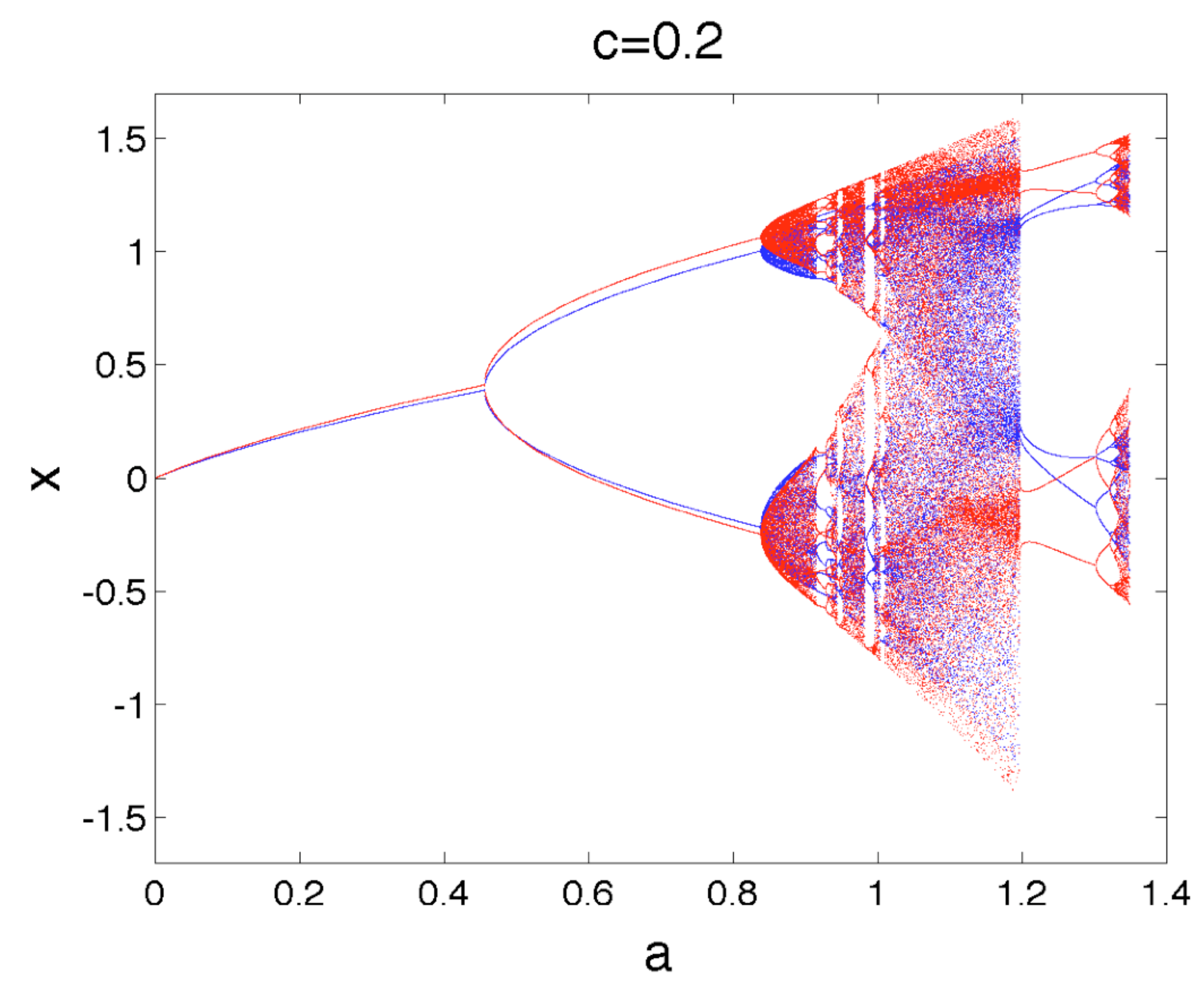}
\includegraphics[width=.45\textwidth]{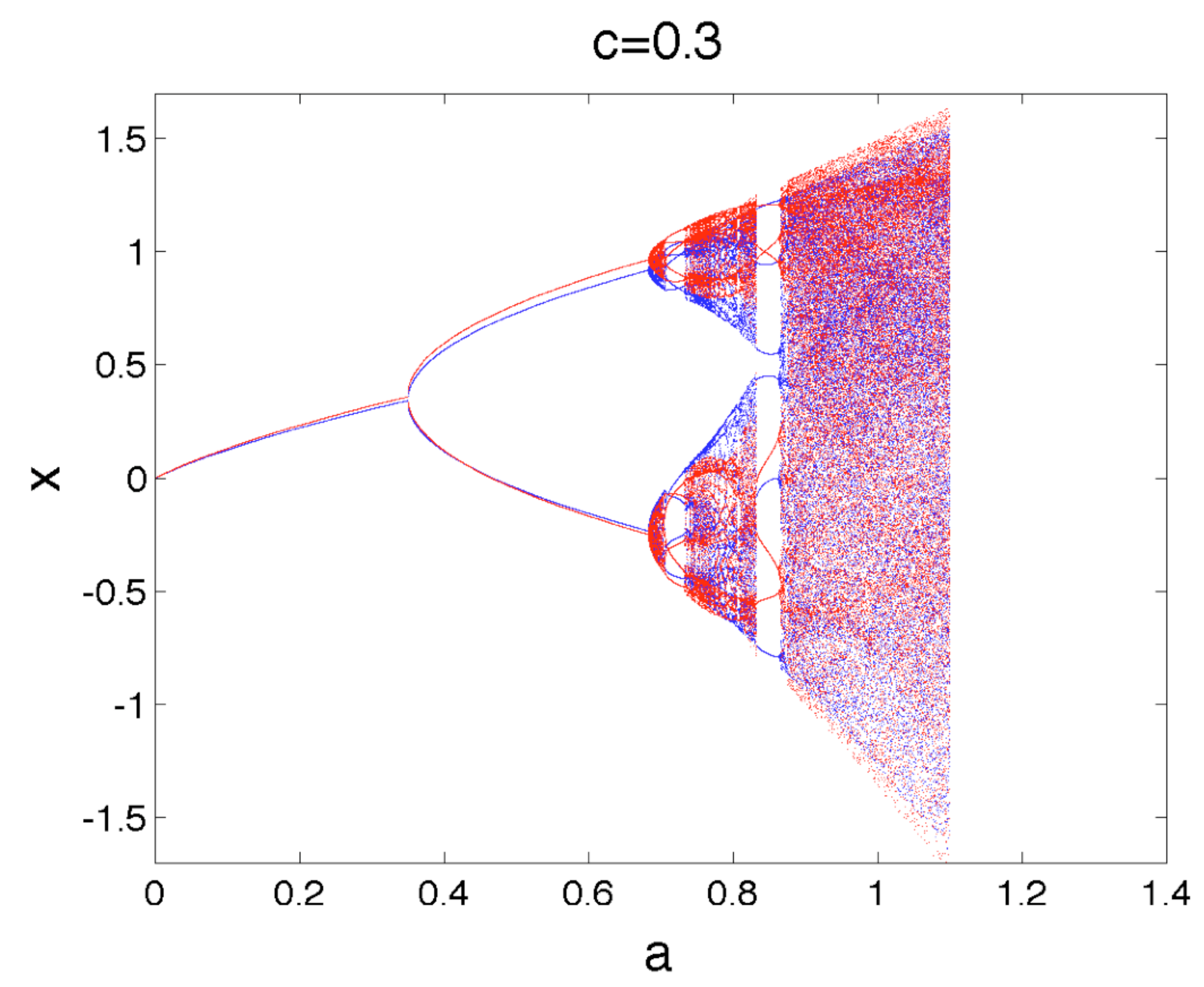}
\caption{\label{fig:coupling}
 Coupled quadratic maps: The attracting set for a pair of coupled
 quadratic maps as a coordinated cascade
 forms. In Section~\ref{sec:coupled}, we prove that infinitely many simultaneously
 occurring cascades occur, though it may be
 that none of them are attracting. The figures above show the stable set
 as a function of parameter $a$ for the system $F(a,(x,y))=(a - x^2 +
 c y, 1.1 \; a -y^2 + c x)$ for four different values of the
 coupling constant $c$. The values of $x$ are depicted in blue, and of
 $y$ in red. For (a) $c= 0$ (no coupling so the bifurcations of the x and y coordinates need not synchronize), (b) $c=0.1$, (c) $c=0.2$, and (d)
 $c=0.3$. For each $c$ and $a$, we depict asymptotic behavior using $(0,0)$
 as initial condition at $a=0$.
}
\end{center}
\end{figure}

\bigskip
{\bf Corollary~\ref{cor:QuadSys}: $N$-dimensional coupled systems.}  This is an extension of the result of Corollary~\ref{cor:1D-quadratic} to the case of $N$ coupled quadratic maps. We describe the number of cascades. See
Figure~\ref{fig:coupling}.   Such cascades are a collective
phenomenon for high-dimensional coupled systems.

\bigskip

{\bf The literature on general existence of cascades is scant.}
M. Feigenbaum~\cite{feigenbaum:79} developed methods that have been
used to rigorously demonstrate the existence of a cascade for a
parametrized map when it is sufficiently close to quadratic.  Indeed
his goal was to show that when a period-one cascade exists for a
parametrized map that is nearly quadratic, the procession of
period-doubling values has a regular scaling behavior.  In contrast,
our goal is to establish criteria for the existence of infinitely many
cascades.  It seems likely that Feigenbaum's scaling rule holds for
typical cascades of the parametrized maps in this paper. However, our
methods are topological and say nothing about scaling.

In the theory of one-dimensional quadratic maps, the existence of cascades became a folk theorem based on the property of {\bf monotonicity}: Namely, that as the parameter increases, new orbits can appear but no orbits are destroyed.
This monotonicity was originally proved by Douady and Hubbard in the complex analytic setting. See 
\cite{milnor:thurston:88} for a proof. When there is monotonicity, the existence of cascades is quite straightforward. Monotonicity of the quadratic map implies that the only possible periodic-orbit bifurcations are those saddle-node and period-doubling bifurcations in which periodic attractors are created.
Periodic attractors do not persist as attractors as $\lambda$ increases because there are no periodic attractors for $\lambda \ge 2$; in order for a periodic attractor to become unstable as $\lambda$ increases, it must undergo a period-doubling bifurcation, creating a new periodic attractor with double the period. This new periodic attractor must cease to exist before $\lambda = 2$, so the new higher period attractor undergoes a period doubling as well, etc. The attractors must undergo infinitely many period-doubling bifurcations as $\lambda$ increases. Hence a cascade exists.

General one-dimensional systems need not be monotonic, and higher-dimensional systems tend never to be monotonic as shown in~\cite{yorke:antimonotonicity-Annals92}. In these cases there are many more possible bifurcations which both create and destroy orbits. 
See Figure~\ref{fig:boundedcascades}. 
However, our results demonstrate that the existence of cascades is in no way dependent on either dimension one, nor on monotonicity, nor on having attractors.

Yorke and Alligood~\cite{yorke:alligood:83} discuss in detail a case where the cascade of
period doublings involves attractors, and so they restrict attention to
the case of systems where trajectories are at most one-dimensionally
unstable. We make no such restriction. In higher dimensions, there is no reason for attractors to
be present in cascades, and attractorless cascades do exist.

\bigskip

{\bf The paper proceeds as follows:} In Section~\ref{sec:background},
we classify the set of generic bifurcations of orbits.  In
Section~\ref{sec:orientation}, we develop the orbit index, in order to
investigate the ``index orientation'' on each component.
Section~\ref{sec:abstracttheorem} contains Theorem~\ref{thm:main}.  
In Section~\ref{sec:examples}, we
apply our abstract theorem in Corollaries 1-5 to show that cascades
occur in classes of low-dimensional parametrized maps.  In
Section~\ref{sec:coupled}, we describe classes of $N$-dimensional
parametrized maps that have infinitely many cascades, and we report
the number of cascades of each period.

\newpage

\section{Components in the space of orbits}\label{sec:background}

This section introduces formal versions of the main concepts described
in the introduction, including the space of orbits under the Hausdorff
metric, the set of flip orbits, the set of nonflip orbits, and
cascades viewed as subsets of this space. Throughout this paper we
will assume the following context. To avoid technicalities, we follow
a $C^{\infty}$ approach similar to Milnor's treatment of Sard's
Theorem~\cite{milnor:65}, proving our results for $C^{\infty}$ rather
than for $C^r$ for a specific $r$.

\begin{hypothesis}[The setting]\label{hyp:1}
 Let $F:R \times {\mathfrak M} \to {\mathfrak M}$ be $C^{\infty}$- smooth. We refer to it as a {\bf parametrized map} on ${\mathfrak M}$.
\end{hypothesis}

\begin{definition}[Orbits, flip orbits, and nonflip orbits]\label{dpof}
Write $[x]$ for the orbit of the periodic point $x$. In this paper,
{\bf orbit} always means {\bf periodic orbit}.
By {\bf period} of an orbit or point, we mean its  least period. If $x$ is a periodic point for $F(\lambda, \cdot )$, then we sometimes say $\sigma = (\lambda, x)$
is a periodic point and write $[\sigma]$ or $(\lambda, [x])$ for its orbit.
Let $\sigma =(\lambda, x)$ be a periodic point of period $p$ of a smooth map $G = F(\lambda, \cdot)$.
We refer to the {\bf eigenvalues of $\sigma$ or $[\sigma]$} as shorthand for the eigenvalues of
Jacobian matrix $DG^p(x)$. Of course all the points of an orbit have the same eigenvalues. We say that $[\sigma]$ is {\bf hyperbolic} if none of its eigenvalues have absolute value $1$.
We say it is a {\bf flip orbit} if the number of its eigenvalues (adding multiplicities) less than $- 1$ is odd, and $- 1$ is not an eigenvalue.
We call all other orbits {\bf nonflip orbits}.
Define \[ PO(F)= \left\{ [\sigma]: [\sigma] \mbox{ is an orbit for } F \right\} , \mbox{ and } \]
\[PO_{nonflip}(F) = \left\{ [\sigma] \in PO(F): [\sigma] \mbox{ is a nonflip orbit for } F \right\}.\]
\end{definition}

The distance between two orbits in the space $PO(F)$ or
$PO_{nonflip}(F)$ is defined using the Hausdorff metric.  We say that
two orbits are close in $ {\mathfrak M} $ if every point of each orbit
is close to some point of the other orbit. The periods of the two
orbits need not be the same. This is made precise in the following
definition.

\begin{definition}[Hausdorff metric on sets]
For a compact set $S$ of ${\mathfrak M}$ and $\epsilon > 0$, let $B(\epsilon, S)$ be the closed $\epsilon$ neighborhood of $S$.
Let $S_1$ and $S_2$ be compact subsets of ${\mathfrak M}$.
(We are only interested in the case where these sets are orbits.)
Assume $\epsilon$ is chosen as small as possible such that
$S_1 \subset B(\epsilon,S_2)$ and $S_2 \subset B(\epsilon,S_1)$.
Then the {\bf Hausdorff distance} $dist(S_1, S_2)$ between $S_1$ and $S_2$ is defined to be $\epsilon$.

Let $\sigma_j = (\lambda_j, x_j)$ for $j = 1,2$ be orbits.
We define the distance
between $[\sigma_1]$ and $[\sigma_2]$ to be
\[
\mbox{dist}([\sigma_1], [\sigma_2]) = \mbox{dist}([x_1], [x_2]) + |\lambda_1 - \lambda_2|.
\]
\end{definition}

For example, if $(\lambda, [x(\lambda)])$ is a family of period-$2p$ orbits that bifurcate from the period-$p$ orbit
$(\lambda_*, [x_*])$ due to a period-doubling bifurcation, then
\[
\mbox{dist}((\lambda, [x(\lambda)]), (\lambda_*, [x_*])) \to 0 \mbox{ as } \lambda \to \lambda_*.
\]

\begin{remark} {\em Generally one expects periodic points to be dense
    in a compact chaotic set. Using the Hausdorff metric changes the
    geometry.  Every $R\times {\mathfrak M}$ neighborhood of a saddle
    fixed point $(\lambda,x)$ with a transverse homoclinic
    intersection has infinitely many periodic points $y_m$ (of
    arbitrarily large period). However, since $x$ is a hyperbolic
    saddle point, some points in the orbit of each $y_m$ are far from
    $x$. Thus the orbits $[y_m]$ do not converge to $[x]$ in the
    Hausdorff metric.}
\end{remark}

\begin{definition}[Cascade of period $m$]\label{def:cascade}
The term {\bf component} means a connected component of $PO_{nonflip}(F)$ in the Hausdorff metric.
An {\bf arc} is a set that is homeomorphic to an interval. We call it an {\bf open arc} if the interval is open or {\bf half-open} if that describes the interval.

A {\bf (period-doubling) cascade of period $m$} is a half-open  arc $C$ in
$PO_{nonflip}(F)$ 
with the following properties.
Let $h: [0,1) \to C$ be a homeomorphism.
\begin{enumerate}
\item[(i)] The set of periods of orbits in $C$ is ${m, 2 m, 4 m, 8 m, \cdots}$.
\item[(ii)] The number $m$ is the minimum period of orbits in the component that contains $C$.
\item[(iii)] $C$ has no proper connected subset with properties (i) and (ii).
Note that $h(0)$ will be the only orbit of period $m$ in $C$ and it will be a period-doubling bifurcation orbit.
\end{enumerate}

\bigskip

If a component contains a cascade, we refer to the component as a {\bf cascade component}.
\end{definition}

\begin{remark}
{\em
The cascades we discuss in this paper each lie in a compact subset of $R \times {\mathfrak M}$ (though the cascade's component may be unbounded). 

For generic $F$, we will show that such cascades have the following additional property.
}
\begin{enumerate}
\item[(iv)] If $\{ p_k \}_1^\infty$ is the sequence of periods of the orbits, ordered so that for each $k$ the $k+1$ orbit lies ``between'' (using the ordering induced from [0, 1)) the $k$ orbit and the $k+2$ orbit, then no period will occur more than a finite number of times. That implies 
\[ \lim_{k\to\infty} p_k = \infty.\]
\end{enumerate}
\end{remark}

{\bf Examples of arcs and a cascade.} There is a simple example of a
subset of the orbits in $PO_{nonflip}(F)$ which is homeomorphic to an
interval. Let $F = \lambda - x^2$. The smallest $\lambda$ for which
there is an orbit is $- 1/4$, and that is a saddle-node fixed point $Q
= (- 1/4, - 1/2)$. There is a unique periodic attractor for each
$\lambda \in J = (-1/4, \lambda_{Feig})$ where $\lambda_{Feig}$ is the
end of the first cascade. We often refer to $\lambda_{Feig}$ as the
``Feigenbaum limit parameter.'' The attractors for $\lambda \in J$
constitute a component $C$ of the attractors in $PO_{nonflip}(F)$. For
each $\lambda \in J$ there is a unique attracting orbit $x$ in $C$ and
trivially for each orbit in $C$ there is a unique $\lambda$ in $J$. In fact
the map on $C$ defined by
\[
(\lambda,[x]) \mapsto \lambda
\]
is a homeomorphism. Therefore $C$ is an open arc since it is
homeomorphic to the interval $J$. We shall see below that $C$ is not
maximal.

{\bf The cascade.} There is a subarc $C_1 \subset C$ that is a
cascade. For $\lambda =3/4$ the orbit in $C$ is the period-doubling fixed point $Q_1
= (3/4, 1/2)$. For $\lambda > 3/4$, the orbits in $C$ have period
greater than 1.  Let $C_1 = \{$ the orbits in $C$ for which $\lambda
\ge 3/4\}$.  Then the cascade $C_1$ is the smallest subarc of $C$ for which all
periods $2^k$ occur.

{\bf The component containing $C$.} The arc $C$ is not a maximal
arc. There is another arc, an arc of unstable (derivative $> +1$)
fixed points $(\lambda,y(\lambda))$ defined for $\lambda \in (- 1/4,
+\infty),$ where $y(\lambda) = (-1+ \sqrt{1 + 4 \lambda})/2$.  The two
arcs terminate at the saddle-node fixed point $Q$. Taking the union of
the two arcs plus ${Q}$ yields the component containing $C$. It is 
maximal because on one extreme the period goes to $\infty$ and on the
other extreme, $\lambda \to \infty$. Each point of the arc is a
different orbit. Since the set of $\lambda$ values in the maximal arc
is unbounded, the cascade is an unbounded cascade.

In contrast to this straightforward example for quadratic maps, in
which there is exactly one saddle-node bifurcation, and exactly one
period-doubling bifurcation for each period $2^k$, a cascade is
generally quite complicated. For example, a bounded component that is
a open arc contains two disjoint cascades (one on each end). See
Figure~\ref{fig:boundedcascades}. Quadratic maps have a {\bf
monotonicity} property; if an orbit exists at $\lambda_*$, then it
exists for all $\lambda>\lambda_*$. This leads to cascades which
consist only of forward-directed period doublings of attractors. In
contrast, a cascade for a general map generally does not have such
regular behavior, nor does it consist only of attractors.  In what
follows we will be interested in arcs of orbits with no regard to
whether the orbits are attractors.

\begin{figure}[thb]
\begin{center}
\includegraphics[width=.45\textwidth]{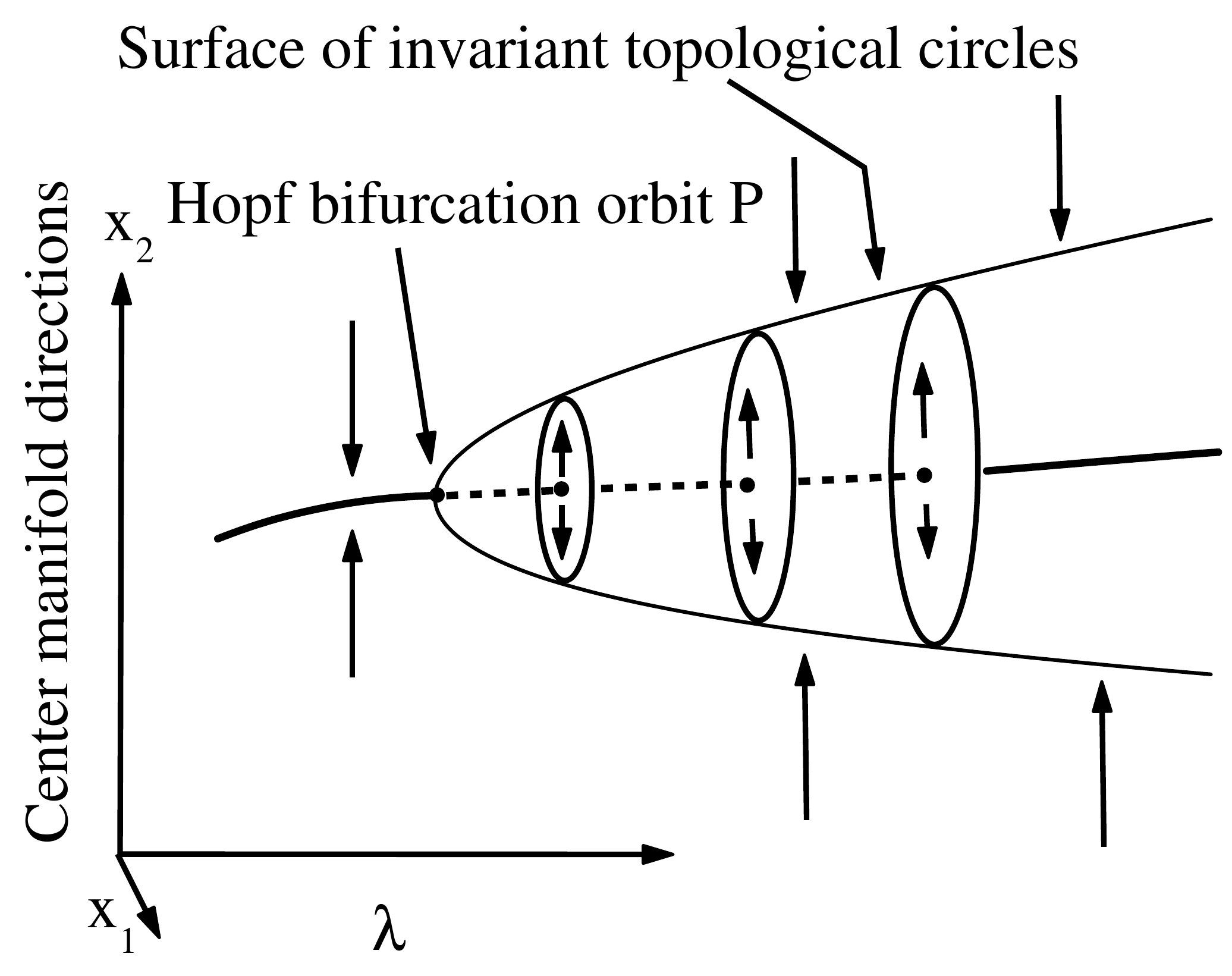}
\caption{\label{fig:hopf1} The Hopf bifurcation: All the local invariant sets of a generic Hopf
 bifurcation lie on a three-dimensional center manifold. There is a
 two-dimensional surface consisting of invariant topological circles
 as well as the arc of periodic points from which it bifurcates. As
 indicated by vertical arrows, in the case shown here the invariant
 surface is attracting (when the dynamics are restricted to the
 center manifold). There also exist generic Hopf bifurcations for which the arrows are all reversed. In this depiction, the
 surface appears as $\lambda$ increases, but it could also occur as
 $\lambda$ decreases.}
\end{center}
\end{figure}

We now explain what is meant by a generic orbit bifurcation and show there is a residual set of $F$ for which every orbit is either hyperbolic or is a generic bifurcation orbit (Proposition~\ref{prop:genfamily} below).

\subsection{Generic Orbit Bifurcations}

\begin{figure}
\begin{center}
\includegraphics[width=.7\textwidth]{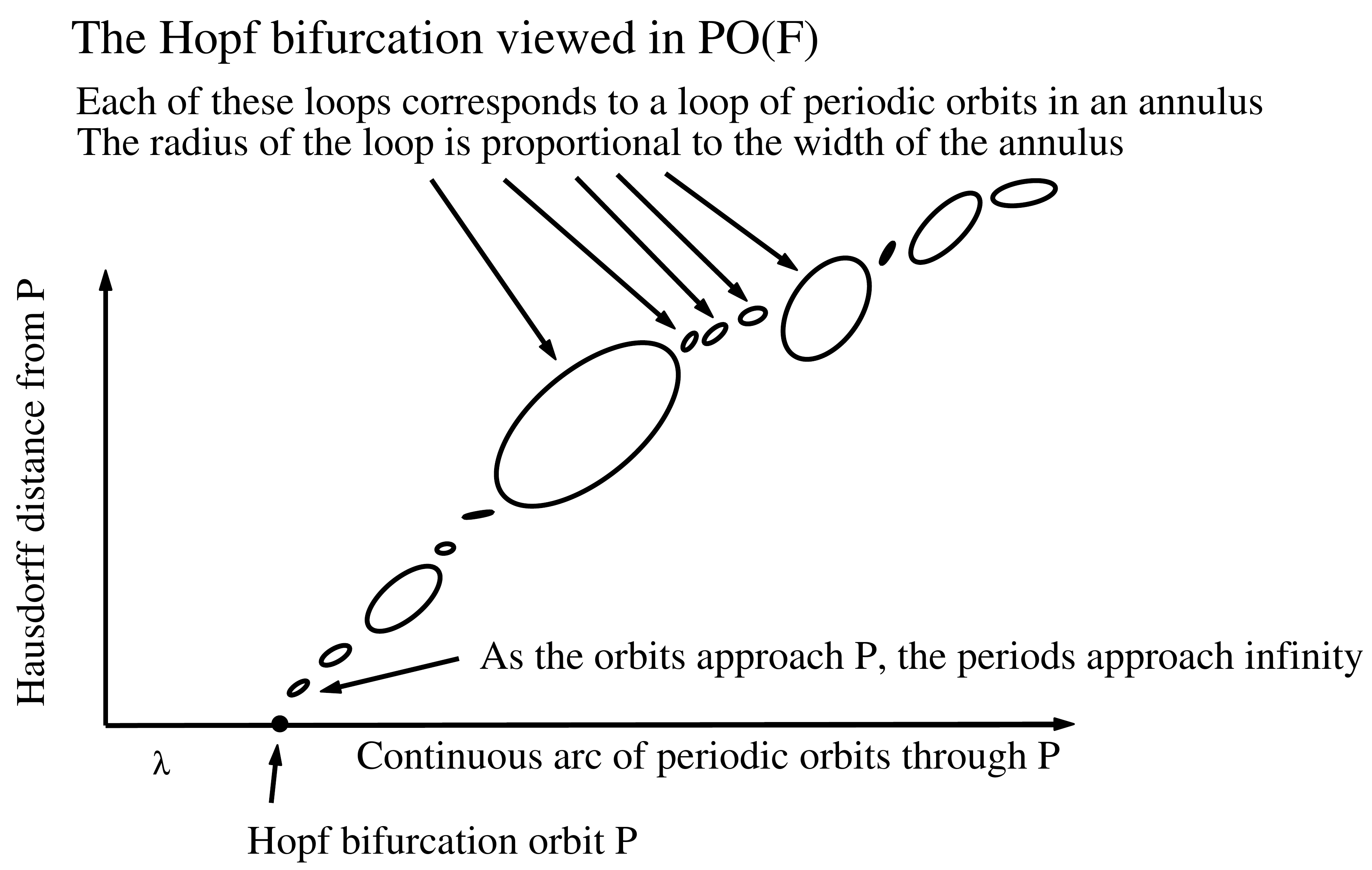}
\caption{\label{fig:hopf4} Orbits near a Hopf bifurcation in $PO(F)$: A neighborhood in $PO(F)$ of a generic Hopf bifurcation orbit $P$ consists of (i) an arc of orbits (the horizontal axis) with the same period as $P$, and (ii) a collection of components, each of which is a loop of orbits, i.e., a simple closed curve. All the orbits in each of these loops have the same period. As the loops approach $P$ in $PO(F)$, the periods go to infinity.
}
\end{center}
\end{figure}

We list the three kinds of generic orbit bifurcations, largely following the treatment in Robinson~\cite{robinson:95}. Cases {\it (i)} and {\it (ii)} are depicted in Figure~\ref{fig:SadNode&PerDoubling}. Figures~\ref{fig:hopf1} and~\ref{fig:hopf2}--\ref{fig:hopf4} depict case {\it (iii)}. 

\begin{definition}[Generic orbit bifurcations]\label{def:genericBifurcation}
 Let $F$ satisfy Hypothesis~\ref{hyp:1}. We say a {\bf bifurcation orbit
 $P$ of $F$ is generic} if it is one of the following three types (as described in detail by Robinson~\cite{robinson:95}):
\begin{enumerate}
\item[\it (i).] A generic saddle-node bifurcation (having eigenvalue $+1$).
\item[\it (ii).] A generic period-doubling bifurcation (having eigenvalue $-1$).
\item[\it (iii).] A generic Hopf bifurcation with complex conjugate eigenvalues which are not roots of unity. 
\end{enumerate}

We say {\bf F is generic} if each non-hyperbolic orbit is one of the above three types.
\end{definition}

Note that the standard Hopf bifurcation theorem permits complex conjugate
eigenvalues to be higher-order roots of unity. However, we have chosen a more stringent
generic bifurcation condition, since a parametrized map with
a bifurcation through complex conjugate pairs that are roots of unity
can be rotated by an arbitrarily small perturbation to a family with a
bifurcation through one complex conjugate pair which are not roots of
unity. Neither 1 nor $-1$ can be perturbed away in this manner, since
these eigenvalues are real, so do not occur in conjugate pairs. Since
the roots of unity are countable, the families with Hopf bifurcations
for which the eigenvalues are not complex roots of unity are generic.
Hence we exclude bifurcations such as period tripling or multiples
other than two.

Further, a generic Hopf bifurcation orbit has a neighborhood in
$PO(F)$ in which the only bifurcation orbits are saddle-node
bifurcation orbits. See Figures~\ref{fig:hopf1}--\ref{fig:hopf4}.

\bigskip
\begin{hypothesis}[Generic bifurcations]\label{hyp:2}
Assume Hypothesis~\ref{hyp:1}. Assume that each orbit of $F$ is either hyperbolic or is a generic bifurcation orbit.
\end{hypothesis}

\begin{proposition}[Generic $F$ constitute a residual set]
\label{prop:genfamily}
There is a residual set $S \subset C^\infty$ of parametrized maps $F$ satisfying Hypothesis~\ref{hyp:1} for which all bifurcation orbits are generic.

This residual set $S$ is $C^1$ dense in the uniform $C^1$ topology;
that is, for each $F \in C^\infty$, there is a sequence $(F_i) \subset S$ such that
$||F-F_i||_{C^1} \to 0$ as $i\to \infty$.
\end{proposition}

The proof of this proposition uses standard transversality arguments. See Palis
and Takens~\cite{palis:takens:87}. Minor changes to their methods give
the irrational rotation number for Hopf bifurcations.

\bigskip

\subsection{Components are one-manifolds for generic $F$}

\begin{definition}[Index Orientation]\label{def:IndexOrientation}

Assume that a periodic orbit $y$ of
 period $p$ of a smooth map $G$ is hyperbolic.
Define the {\bf unstable dimension $dim_u(y)$} 
to be the number
 of real eigenvalues (with multiplicity) having absolute value $>1$.

Let $F$ satisfy Hypothesis~\ref{hyp:2}, and assume that $Q$ is a component that is a one-manifold. (We show below that this is true for generic $F$.) Then we know that there is a homeomorphism $h:X \to Q$ where 
$X$ is either the interval $(-1, +1)$ or a circle, which we will write
as $[-1,1]/ \{ -1,1 \}$. 
%%This is substantially revised.
For each $s \in X$, let $h_\lambda(s)$ denote the projection of $h(s)$ to the corresponding parameter value. Thus $h_\lambda$ is a map from $R$ to $R$, and  $X$ and $\lambda$ both inherit an orientation from the real numbers. Therefore we can describe $h$ as increasing or decreasing at $s$ whenever $h_{\lambda}$  is  increasing or decreasing at $s$.

We say the homeomorphism 
$h$ is an {\bf index orientation} for $Q$ if 
whenever $h(s)$ is a hyperbolic orbit,
$h_{\lambda}(s)$ is locally strictly increasing when $dim_u(h(s))$ is odd
and is locally strictly decreasing when $dim_u(h(s))$ is even.
\end{definition}

We now state the main result of this section. 

\begin{theorem}[Components 
are oriented one-manifolds for generic $F$]
\label{thm:OneManifold}
Consider all $F$ as in Hypothesis~\ref{hyp:1}.
There is a residual set of such $F$ for which each component (of $PO_{nonflip}(F)$)
\begin{enumerate}
\item[($M_1$)] 
is a one-manifold, i.e. is either a simple closed curve or is homeomorphic to an open interval, and
\item[($M_2$)]
it has an index orientation.
\end{enumerate}
\end{theorem}

The proof of $(M_1)$ consists of two parts.  First, we show that
in a neighborhood of any hyperbolic orbit, its component is an
arc. Secondly we show the same for each non-hyperbolic orbit. Hence
each point in a component has a neighborhood in the component that is
an arc, so the component is a one-manifold.

{\bf The neighborhood of a hyperbolic orbit.} If $\sigma_0 =
(\lambda_0, x_0)$ is a period-$p$ hyperbolic point, the implicit
function theorem implies that it has a smooth unique {\bf continuation
  curve} of period-$p$ points $\sigma(\lambda) = (\lambda,
x(\lambda))$, defined for $\lambda$ in some neighborhood of
$\lambda_0$ with $\sigma(\lambda_0) = \sigma_0$. Furthermore
$[\sigma_0]$ has a neighborhood in $PO(F)$ in which there are no
orbits other than those of the continuation.  Hence each hyperbolic
orbit has a neighborhood in $PO(F)$ which is an arc. 

{\bf The neighborhoods of a non-hyperbolic orbit.}  This part of the
proof relies on the fact that $F$ is generic as defined in Definition~\ref{def:genericBifurcation}.

The proof of Part ($M_1$) will be complete with Proposition~\ref{prop:hopf}, which shows that each nonflip generic bifurcation orbit $y$ has a neighborhood that -- when intersected with the component it is in -- is an arc that passes through $y$.

The proof of Part ($M_2$) of Theorem~\ref{thm:OneManifold} is complete when we finish proving Proposition~\ref{prop:consistentOrientation} in Section~\ref{sec:orientation}.
There we show that for each homeomorphism $h$ from $X$ to a component is either an index orientation or $h^*$ is where $h^*(s) = h(-s)$. Recall that $-s \in X$ whenever $s \in X$.

\begin{figure}
\begin{center}
\includegraphics[width=.5\textwidth]{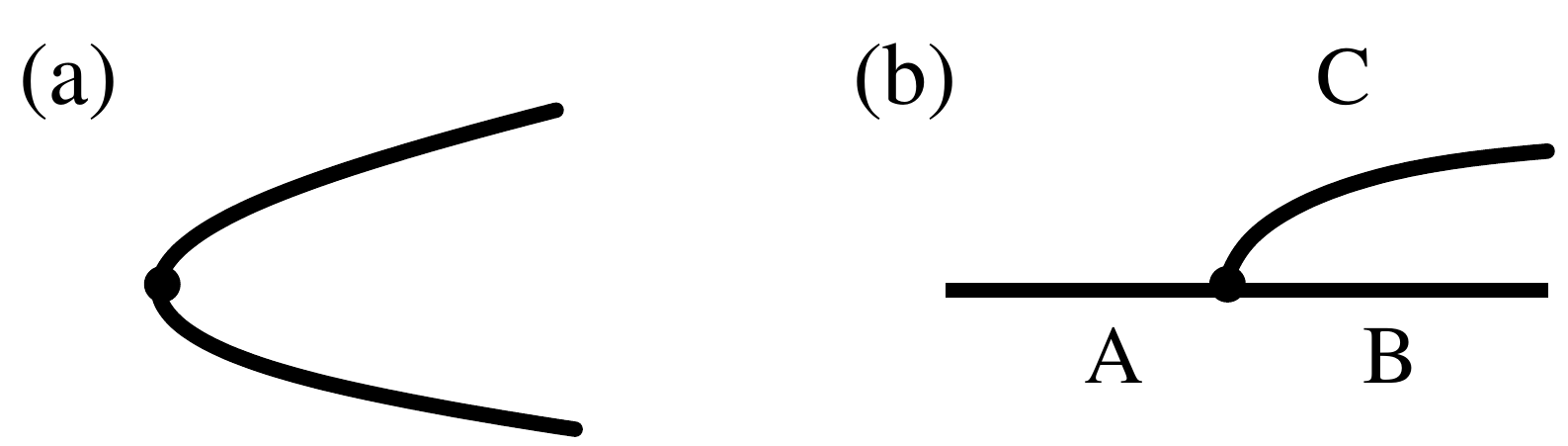}
\caption{\label{fig:SadNode&PerDoubling} 
A depiction of a sufficiently small neighborhood of (a) the saddle-node bifurcation and (b) the period-doubling bifurcation. The horizontal axis is the parameter, but it can be either increasing or decreasing in this figure. The vertical axis is the space $PO(F)$, so each orbit is depicted as a single point. Near the bifurcation point in (a) all orbits are either flip orbits or all orbits are nonflip orbits. 
Near the bifurcation point in (b), exactly one of segments A and B consists of flip orbits, and the other consists of nonflip orbits. The period-doubled segment C always consists of nonflip orbits. Hence exactly two of the segments of a period-doubling bifurcation are contained in $PO_{nonflip}(F)$.
}
\end{center}
\end{figure}

\begin{proposition}[The neighborhood of a generic bifurcation orbit]
\label{prop:hopf}
Assume $F$ satisfies Hypothesis~\ref{hyp:2}. Assume $P=(\lambda_0,[x_0])$ is a generic bifurcation orbit 
in $PO_{nonflip}(F)$ 
and let $C$ be its component.
Then $P$ has a neighborhood in $C$
that is an open arc in which it is the only bifurcation orbit.
\end{proposition}

\begin{figure}
\begin{center}
\includegraphics[width=0.55\textwidth]{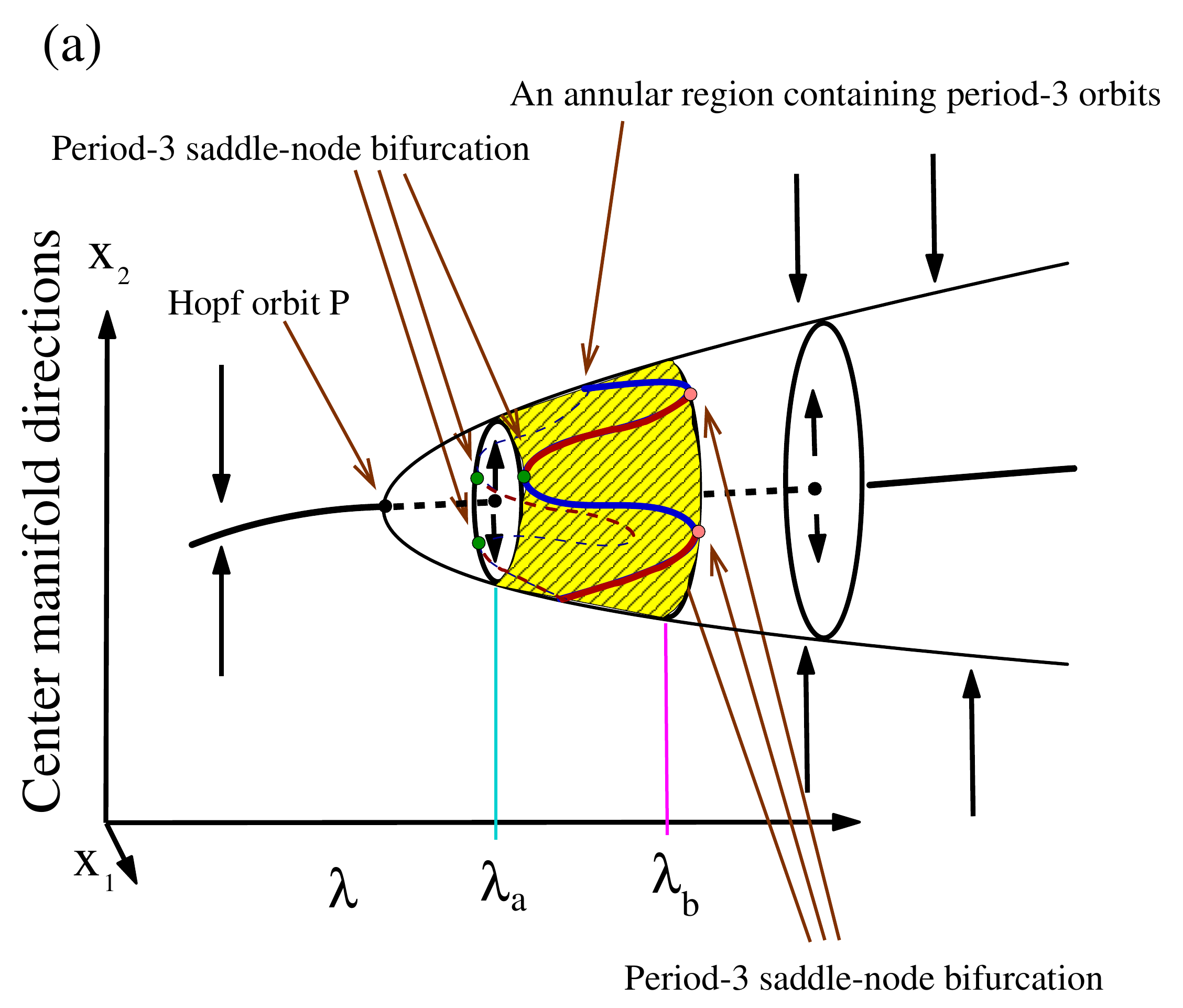}
\includegraphics[width=0.6\textwidth]{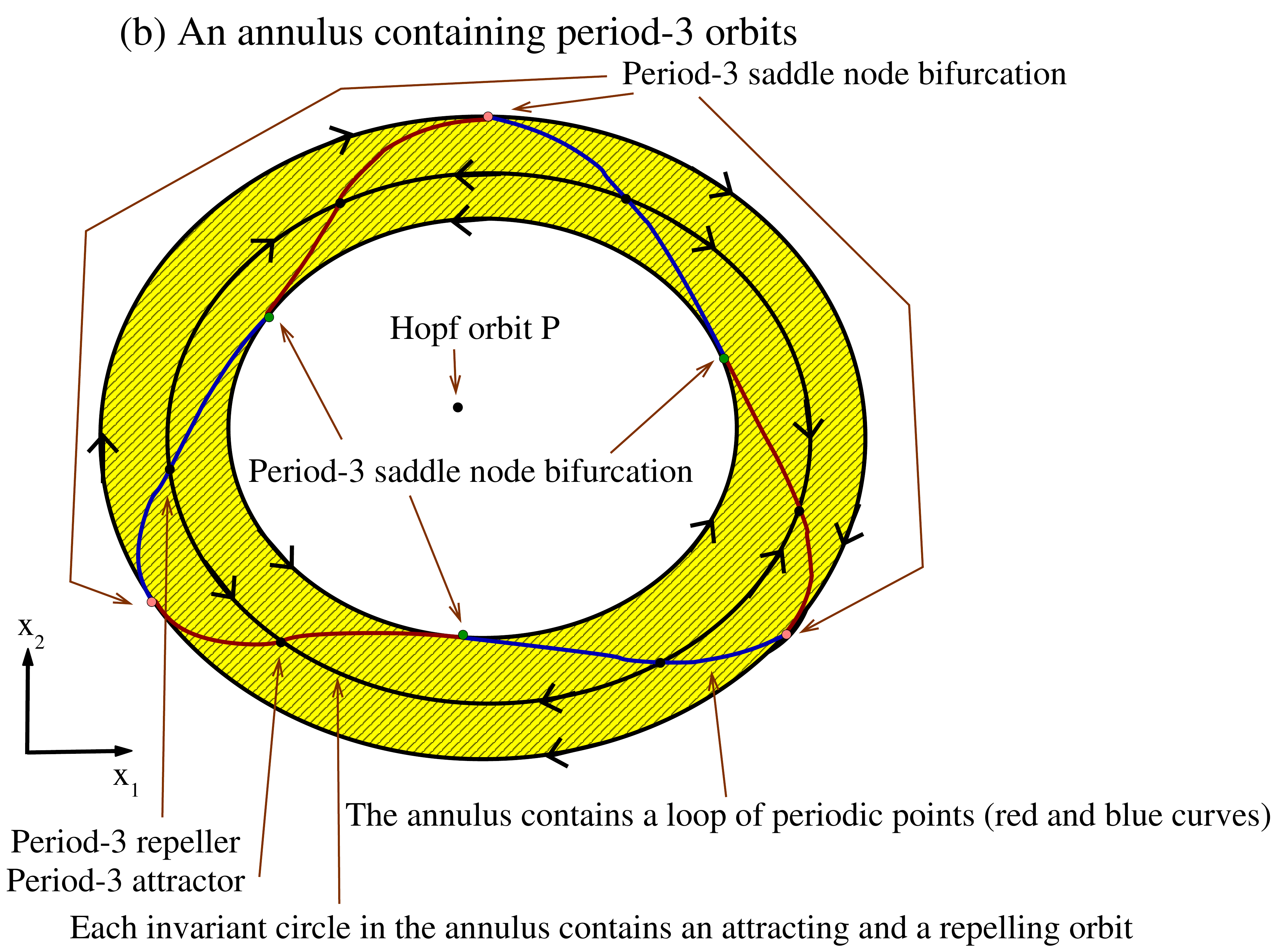}
\includegraphics[width=0.39\textwidth]{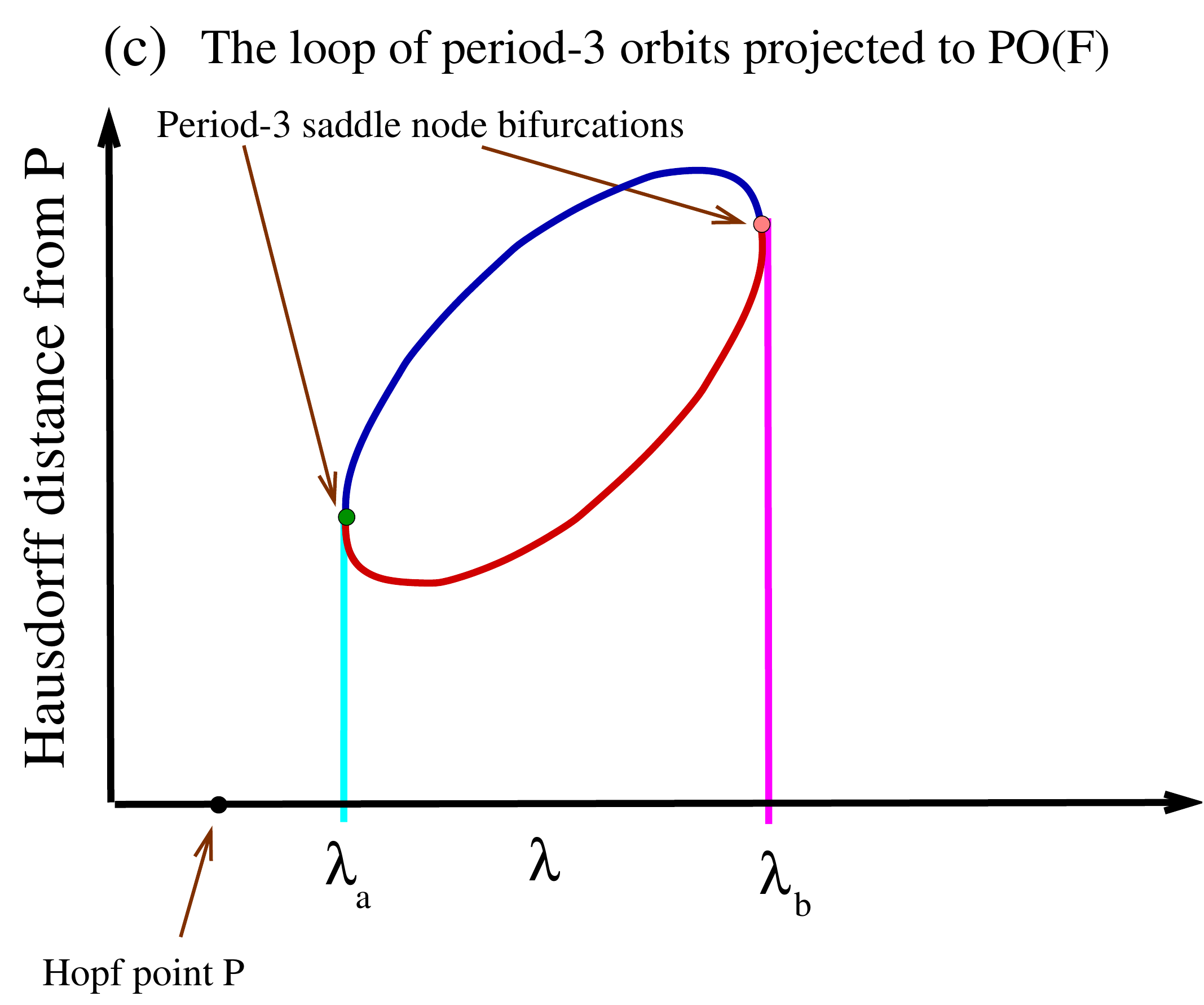}
\caption{\label{fig:hopf2} Orbits near a Hopf bifurcation: (a) Within the two-dimensional surface of
 invariant circles near a generic Hopf bifurcation, the topological invariant circles containing
 orbits of a fixed period form annular regions. (b) Each annular
 region projects to an annulus when projected to the plane of spatial
 directions of the center manifold. The annulus consists of
 invariant topological circles, and each of those circles has an
 attracting period-$k$ orbit (attracting in the circle but not in the
 annulus) and a repelling period-$k$ orbit - except for the inner and
 outer boundary circles of the annulus. The boundary circles contain
 bifurcating period-$k$ saddle-node orbits. (c) The loop of periodic points is a $k$-fold cover of the
 corresponding loop of period-$k$ orbits in $PO(F)$. The component
 of the period-$k$ orbits is a topological circle (that is, a loop) in $PO(F)$. }
\end{center}
\end{figure}

\begin{figure}
\begin{center}
\includegraphics[width=.6\textwidth]{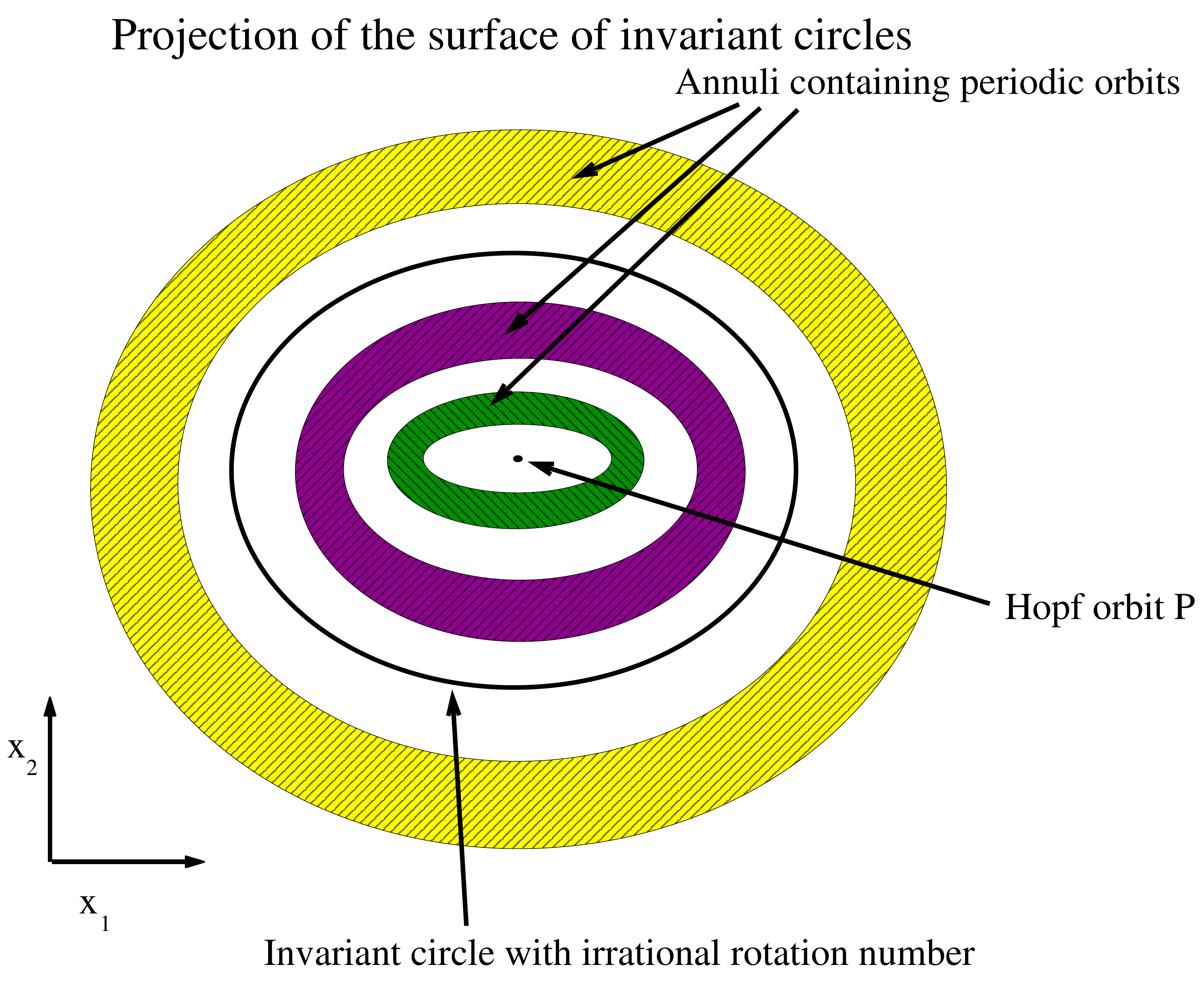}
\caption{\label{fig:hopf3}
 Orbits near a Hopf bifurcation lying on the two-dimensional surface in the three-dimensional center manifold: Here we project the parabolic region shown in
Figure~\ref{fig:hopf2} onto a plane. Generically, near a Hopf bifurcation orbit, there are infinitely many annular regions of orbits, each of a fixed
 period, separated by invariant circles with irrational rotation
 number. }
\end{center}
\end{figure}

\begin{proof}
We consider each of the three generic orbit bifurcations individually.

{\bf Case (i)}. See Figure~\ref{fig:SadNode&PerDoubling}a. Locally a generic saddle-node orbit $P$ is in an open arc of otherwise hyperbolic orbits. If $P$ is nonflip, then so are the nearby orbits. No other orbits are nearby, so the proposition's assertion is true for these orbits. 

{\bf Case (ii)}. See Figure~\ref{fig:SadNode&PerDoubling}b.
 Locally a generic period-doubling orbit $P$ of period $p$ consists of an arc of period-$p$ orbits passing through $P$, plus an arc of period $2p$ that terminates at $P$. The latter is always nonflip near $P$. The period-$p$ branch has an eigenvalue that passes through $-1$ at $P$. Hence (locally) on one side of $P$ the arc consists of nonflip orbits and flip orbits on the other side. Hence there are two arcs in $PO_{nonflip}(F)$ that terminate at $P$. Hence a neighborhood of $P$ in $PO_{nonflip}(F)$ is an arc.

{\bf Case (iii)}. The generic Hopf bifurcation is much more complicated. It is the only generic local bifurcation for which it is possible to have orbits of unbounded periods limiting to $P$ in the Hausdorff metric.

 The center manifold theorem guarantees that for any $r$ in a
 sufficiently small neighborhood of $P$, there is a three-dimensional $C^r$ center
 manifold for $P$ in $R \times {\mathfrak M}$. The Hopf Bifurcation Theorem
 guarantees that within this center manifold, there is an invariant topological
 paraboloid which is either attracting, as depicted in Figure~\ref{fig:hopf1}, or repelling. Within this paraboloid, for each $\lambda$ value there is
an associated rotation number $\omega_\lambda$ of the invariant circle at parameter value
 $\lambda$, as shown in Figures~\ref{fig:hopf2} and~\ref{fig:hopf3}.

 By the genericity assumption on Hopf bifurcation points, the limit $\omega_{\lambda_0}$ as
 $\lambda \to \lambda_0$ of the rotation numbers $\omega_\lambda$ is
 irrational. If $\omega_\lambda$ is constant, then there are no
 local orbits, so $P$ would not be a bifurcation orbit. If
 $\omega_\lambda$ is non-constant, then it varies through an interval
 which contains both irrational and rational values. Let $Y$ be
 the local set of orbits in $PO(F)$ other than $\left\{
 (\lambda,[x_\lambda]) \right\}$, the continuation of $P$. All
 points in $Y$ are on the invariant paraboloid. Therefore
 $(\lambda,[y]) \in Y$ is only possible for $\lambda$ values such
 that $\omega_\lambda$ is rational. Therefore points in $Y$ are disconnected.

 Specifically, take any point in $Y$. It is not in the same component
 as $\left\{ (\lambda,[x_\lambda]) \right\}$ since the paraboloid can be separated into two components at
 every parameter value for which the rotation is irrational.
 Therefore the curve $\left\{ (\lambda,[x_\lambda])
 \right\}$ is isolated in its component of $PO(F)$.
Furthermore, since a Hopf bifurcation changes the number of eigenvalues outside the unit circle by two, the curve $\{\lambda,[x_\lambda]\}$ is either entirely flip orbits or is entirely contained in $PO_{nonflip}(F)$.
\end{proof}

\begin{remark}\label{lemmaub} {\em Under Hypothesis~\ref{hyp:2}, the period of the orbits in a component $C$ is locally constant near hyperbolic orbits and near saddle-node and Hopf bifurcations. The period
can only change at period-doubling bifurcations, in which case it changes by a factor of two. Hence an arc $C$ in $PO_{nonflip}(F)$ is a cascade only if the sequence of periods $\left\{ p_k \right\}$ of the non-hyperbolic orbits in $C$ (such that orbit $k+1$ is between orbit $k$ and orbit $k+2$) limits to infinity. This occurs if the sequence is infinite and no period occurs more than a finite number of times.
}
\end{remark}

\subsection{Bounded arcs and cascades}

For a generic map $F$, let $A$ be a component that is an open arc. 
Then there is a homeomorphism
$h: (-1,1) \to A$. Let $m$ denote the minimum period of the orbits in
$A$. Without loss of generality, we can assume that $h(0)$ is a hyperbolic orbit
with period $m$. Write $A^- = h((-1,0))$ and $A^+ = h((0,+1))$. We
say that a set of orbits in $PO_{nonflip}(F)$ is {\bf bounded} if the union
of its orbits lies in a compact subset of $R\times {\mathfrak M}$. We say an
open arc $A$ has a {\bf bounded end} if either $A^-$ or $A^+$ is
bounded. While we are splitting $A$ at a rather particular orbit, the
property of boundedness of an end is independent of where the arc $A$
is split. Note that by our assumption, the point $h(0)$ is not in a cascade, since a cascade contains only one orbit of smallest period, meaning that it is not hyperbolic. We refer to $A^-$ and  $A^+$ as the {\bf ends of the component}.

\begin{proposition}[Bounded Cascades]\label{prop:bddArcs}
  Assume Hypothesis~\ref{hyp:2}. If a component $A$ is an open arc, and one of its ends is bounded, then that end contains a
  cascade. If the entire component $A$ is
  bounded, then $A$ contains two cascades, and these are disjoint.
\end{proposition}

Write $Per(u)$ for the period of an orbit $u \in PO(F)$. We adopt the
above notation for $m, h, A, A^-$ and $A^+$. For
brevity, we give the proof for the case where $A$ is bounded, since 
the case of one end being bounded uses the exact same method of proof.

\begin{proof} 
 Write $J = (-1,1)$. Let $(t_j)_1^\infty \subset J$ be a sequence
 that converges to either $+1$ or $-1$. Since $A$ is a 
 component, it is closed in $PO_{nonflip}(F)$, so the sequence has no limit points.
 Let $m_j$ denote $Per(h(t_j))$ for each $j$. If a subsequence of
 $(m_j)$ were bounded, $h(t_j)$ would have a limit point, which it
 does not, so $\lim_{j\to \infty} m_j = \infty$. Let $m = \min_j \{
 m_j \}$. Let $S = {t \in J : Per(h(t)) = m}$. Note that S
 is compact. Let $t_{\sup} = \sup S$ and
 $t_{\inf} = \inf S$. 
Note that $-1 < t_{\inf} < 0 < t_{\sup} < 1$. 
Write $J_1 := (-1,
 t_{\inf}]$ and $J_2 := [t_{\sup},+1)$. We claim $A_1 = h(J_1)$ and
 $A_2 = h(J_2)$ are cascades. Note that they are disjoint and each
 is homeomorphic to a half-open interval under $h$.
 Notice that if the $Per(h(t))$ changes discontinuously at $t$, then
 $h(t)$ is a period-doubling orbit and the change is precisely by a
 factor of 2 ({\it cf.} Proposition~\ref{prop:hopf}). Hence since
 the periods are unbounded in $A_1$ and $A_2$, the periods of orbits
 in each must be $\{ 2^k m: k = 0, 1, 2, \cdots \}$, as required by
 item ({\it i}) in Definition~\ref{def:cascade}. Since $m$ is the smallest period
 in $A$, item ({\it ii}) is also satisfied. Item ({\it iii}) is
 satisfied since by our choice of $t_{\inf}$ and $t_{\sup}$, $A_1$ and
 $A_2$ both have only one orbit with period $m$. The two cascades
 are disjoint because $J_1$ and $J_2$ are disjoint.

We note that $Per(h(t)) \to \infty$ as $|t| \to 1$, as mentioned in item {\it (iv)}.
\end{proof} 

The obvious problem with this result is that we are not told how to
demonstrate that components that are bounded arcs. The key to solving that problem lies in the next section which
provides a natural orientation for each component using the index orientation.
 
Note that one interpretation of the fact that $A$ has no limit points is 
the following property:

{\bf Isolation of generic bifurcation orbits of period $\le p$. }
Assume Hypothesis~\ref{hyp:2}. For each period $p$, $F$ has at most a
finite number of non-hyperbolic orbits of period $p$ in each bounded
region of $R \times {\mathfrak M}$. To prove this, note that if there were an
infinite number of bifurcation orbits of period $\le p$ in a bounded
set, then there would be an accumulation point of these bifurcation
orbits in $PO_{nonflip}(F)$, which would have to be
non-hyperbolic. However, this cannot occur, since generic orbit
bifurcations are isolated from bifurcation orbits of bounded period.

\section{An orientation for components\label{sec:orientation}}

An arc in $PO_{nonflip}(F)$ has two  orientations. This
section establishes that one of these two orientations is
consistent with a specific topological invariant called the orbit
index. We establish the behavior of this index near each generic
bifurcation. From this, we are able to conclude that cascades
occur, detailed in the main theorem in Section~\ref{sec:abstracttheorem}.

The following concept of the {\bf orbit
 index} was developed in~\cite{malletparet:yorke:82}, where it is defined for all isolated orbits (for flows).

\begin{definition}[Orbit Index]
Assume that an orbit $y$ of
 period $p$ of a smooth map $G$ is hyperbolic.
Based on the eigenvalues of $y$, we define
 \begin{eqnarray*}
 \sigma^+=\sigma^+(y)
&=& \mbox{ the number
 of real eigenvalues (with multiplicity) in } (1,\infty).\\
 \sigma^- = \sigma^-(y) &=& \mbox{ the number of real eigenvalues (with multiplicity) in } (-\infty, -1).
\end{eqnarray*}

The {\bf fixed point index} of $y$ is defined as $ind(y) = (-1)^{\sigma^+}$.

From the definition of fixed point index, it follows that 

\begin{eqnarray*} 
ind(x, G^{pm}) &=& (-1)^{\sigma^+} \mbox{ for } m \mbox{ odd,}\\
 &=& (-1)^{\sigma^+ +\sigma^-} \mbox{ for } m \mbox{ even.}
\end{eqnarray*}
 
Since $\sigma^+$ and $\sigma^-$ are the same for each point of an orbit, we can define the
{\bf orbit index} of a hyperbolic orbit

\begin{eqnarray}   \label{e:1*}
\phi([x])=
\left\{
\begin{array}{cl}
(-1)^{\sigma^+} & \mbox{if } \sigma^- \mbox{ is even,}  \\
0 & \mbox{if } \sigma^- \mbox{ is odd.}
\end{array}
\right.
\end{eqnarray}
\end{definition}
Hence if $[x]$ is a nonflip hyperbolic orbit,
\[	\phi([x]) = ind([x]) \]

Note that a hyperbolic orbit is a flip orbit if and only if its orbit index is zero. Thus for every hyperbolic orbit $[x]$ that is nonflip, i.e., $[x] \in PO_{nonflip}(F)$,
$\phi([x])$ is $\pm 1$ (never zero).

The following proposition is a stronger version of part $(M_2)$ of Theorem~\ref{thm:OneManifold} and is used to prove $(M_2)$. 
It states that each component has an
index orientation.

\begin{proposition}[Each component has an index orientation]
\label{prop:consistentOrientation}
Let $F$ satisfy Hypothesis~\ref{hyp:2}, and let $Q$ be a component.
Let $\psi:X \to Q$ be a homeomorphism where $X$ is the circle or interval in Definition~\ref{def:IndexOrientation}.
Define the homeomorphism $\psi^*:X \to Q$ by $\psi^*(s)=\psi(-s)$ for all $s \in X$.
Then either $\psi$ or $\psi^*$ is an index orientation.
\end{proposition}

\begin{proof}
On each nonflip orbit $y$ (i.e., where $\sigma^-(y)$ is even), 
$ \phi(y)= (-1)^{\sigma^+(y)}$. Hence on hyperbolic nonflip orbits, 
$ \phi(y) = (-1)^{dim_u(y)}$. Hence on a hyperbolic nonflip orbit $y$,

\[dim_u(y) \mbox{ is odd if and only if } \phi(y) = -1. \]
Hence in the definition of ``index orientation'' we will substitute ``$ \phi(h(s)) = -1$ (or $+1$)'' for ``$dim_u(h(s))$ is odd (or even, respectively).''

Let $\psi, \psi^*: X \to Q$ be as in the statement of the proposition. 
The component $Q$ consists of pairs of segments of hyperbolic orbits connected at generic bifurcation orbits. Pick $s \in X$ such that $P=\psi(s)=\psi^*(-s)$ is hyperbolic. Since both the direction of the arc and the orbit index are fixed on a hyperbolic segment, either $\psi$ or $\psi^*$ is an index orientation on the segment of hyperbolic orbits containing $P$. Assume without loss of generality that this occurs for $\psi$.
We show below that at each type of orbit bifurcation, there is a ``consistent'' index orientation, as depicted in Figure~\ref{fig:sadd}. 
That is, two consecutive segments of hyperbolic orbits (separated by a bifurcation orbit) have the same index orientation, $\psi$ or $\psi^*$.
One of the two segments always leads toward the
bifurcation orbit and one segment leads away as $s$ increases.
Thus continuing by induction through all of its hyperbolic segments, $\psi$ is consistent with $\phi$ on all hyperbolic segments. We can then conclude that $\psi$ is an index orientation. 

We now must show only that at each type of generic orbit bifurcation, there is a consistent index orientation. 
Here and subsequently, since the indices are fixed on hyperbolic segments, we denote the indices of a segment when we mean the indices of any orbit on that segment.

{\bf Part (i): Saddle-node bifurcations.}
Let $y_0 = (\lambda_0,
[x_0])$ be a generic period-$p$ saddle-node bifurcation orbit. In some neighborhood of $y_0$ such that on one side there will be two
segments of hyperbolic period-$p$ orbits $y_a$ and $y_b$, and on the other side there are no orbits, as depicted in Figure~\ref{fig:SadNode&PerDoubling}a. The $\sigma^+$ values of $y_a$ and $y_b$ will differ
by $1$, since along the arc an eigenvalue passes through $1$, but their
$\sigma^-$ values will be equal. Hence $y_a$ and $y_b$ have opposite fixed point
index $ind(y_a, G^p) = - ind(y_b, G^p)$, both possibly 0. Also
$ind(y_a, G^{2p}) = - ind(y_b, G^{2p})$.  Hence from Eqn.~\ref{e:1*},
\begin{equation}\label{eqn:sumsn} \phi(y_a) + \phi(y_b) = 0. \end{equation}
Thus there is a consistent index orientation at this bifurcation.

\bigskip

{\bf Part (ii): Period-doubling bifurcations.} Let $y_0 = (\lambda_0, [x_0])$ be a generic period-$p$
period-doubling bifurcation orbit. In some
neighborhood of $y_0$, on one side there
will be a segment of hyperbolic period-$2p$ orbits. We denote the orbits by $y_c$, as in Figure~\ref{fig:SadNode&PerDoubling}b. For
$\lambda$ close to $\lambda_0$, $D_xF^{2p}$ is approximately $(D_xF^p)^2$ and so it has
no real eigenvalues less than $-1$, and in particular, $\phi(y_c) \ne 0$. Hence
\begin{equation}\label{eqn:2*}
\phi(y_c) = ind(y_0, G^{2p}).
\end{equation}

On the same side there must be a segment of hyperbolic period-$p$ orbits, which we denote $y_b$.
We write $y_a$ for the segment of period-$p$ orbits on the other side from $y_c$. The
invariance of the total fixed point index at a bifurcation yields
\begin{eqnarray*} 
ind(y_a, G^{p}) &=& ind(y_b, G^{p}), \\
ind(y_a, G^{2p}) &=& ind(y_b, G^{2p}) + 2 \; ind(y_c, G^{2p}).
\end{eqnarray*}

We substitute $\phi$ for $ind$
using Eqn.~\ref{eqn:2*}, take the average of the left sides of the above
equations, and set that equal to the average of the right sides.
This yields
\begin{equation}\label{eqn:sumpd} \phi(y_a) = \phi(y_b) + \phi(y_c). \end{equation}
Since each of these has values in $\left\{ -1,0, +1 \right\}$, and $\phi(y_c)$ is not zero, there are two cases:
\begin{eqnarray*} 
\phi(y_a) &=& 0 \mbox{ and } \phi(y_b) = -\phi(y_c) \mbox{ or}\\
\phi(y_b) &=& 0 \mbox{ and } \phi(y_a) = +\phi(y_c).
\end{eqnarray*}
 
Hence there are two segments on which $\phi$ is nonzero. 
If both nonzero segments
are on the same side of the bifurcation point, they have opposite
orientation. If the segments are on opposite sides of the bifurcation point, they have the same
orientation. 
Thus there is a consistent index orientation at this bifurcation orbit.

{\bf Part (iii): Hopf bifurcations.} If an arc of orbits has a Hopf
bifurcation, $\sigma^+$ and $\sigma^-$ are the same on the two sides
of the bifurcation so $\phi$ does not change. If a pair of complex
values become real as $\lambda$ is varied, $\sigma^+$ or $\sigma^-$
can change by $+2$ or $-2$, which has no effect on 
$\phi$ and $ind$.
Thus there is a consistent index orientation at this bifurcation.

This completes the proof of the proposition.
\end{proof}

\begin{figure}
\begin{center}
\includegraphics[width=.6\textwidth]{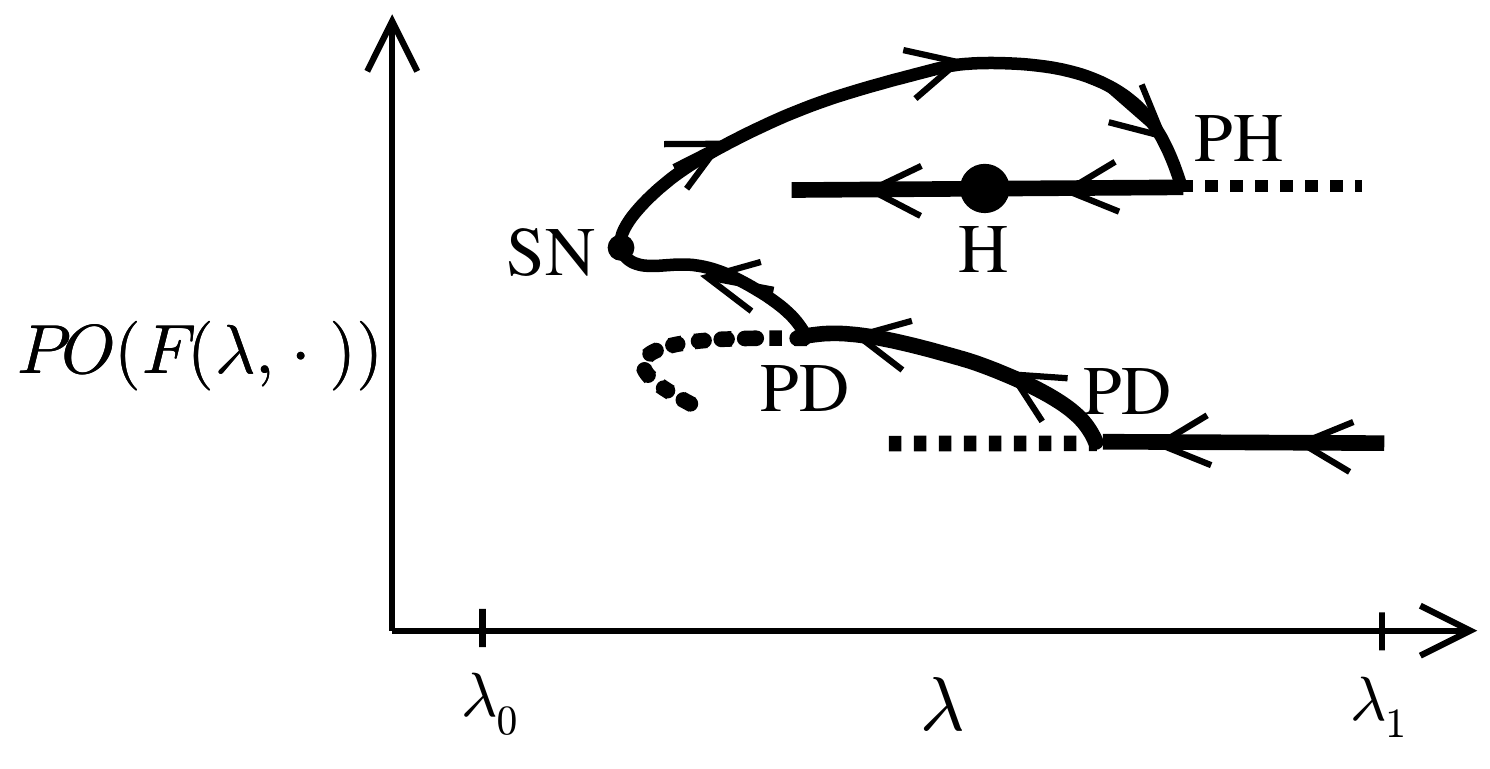}
\caption{\label{fig:snake} Part of oriented component $Q$ with period-doubling
 (PD), period-halving (PH), saddle-node (SN), and Hopf (H)
 bifurcations. If the homeomorphism $h: X \to PO_{nonflip}(F)$ is an index orientation, then as $s \in X$ increases, the $\lambda$ coordinate of $h(s)$ increases or decreases and the arrows on the segments indicate which. Left-pointing arrows correspond to orbit index $\phi = -1$, and right-pointing to $+1$. One of the two adjacent segments always leads toward the intervening
bifurcation orbit and one segment leads away as $s$ increases.
The dotted lines indicate flip orbits (which are not in $Q$).
}
\end{center}
\end{figure}

Figure~\ref{fig:snake} depicts a typical oriented arc in $PO_{nonflip}(F)$ as described in the above result.
A proof similar to the proof of the above proposition also shows that the orbit index is a bifurcation invariant for generic bifurcations.  See \cite{malletparet:yorke:82}. Figure~\ref{fig:sadd} depicts all generic bifurcations.

\begin{figure}
\begin{center}
\includegraphics[width=\textwidth]{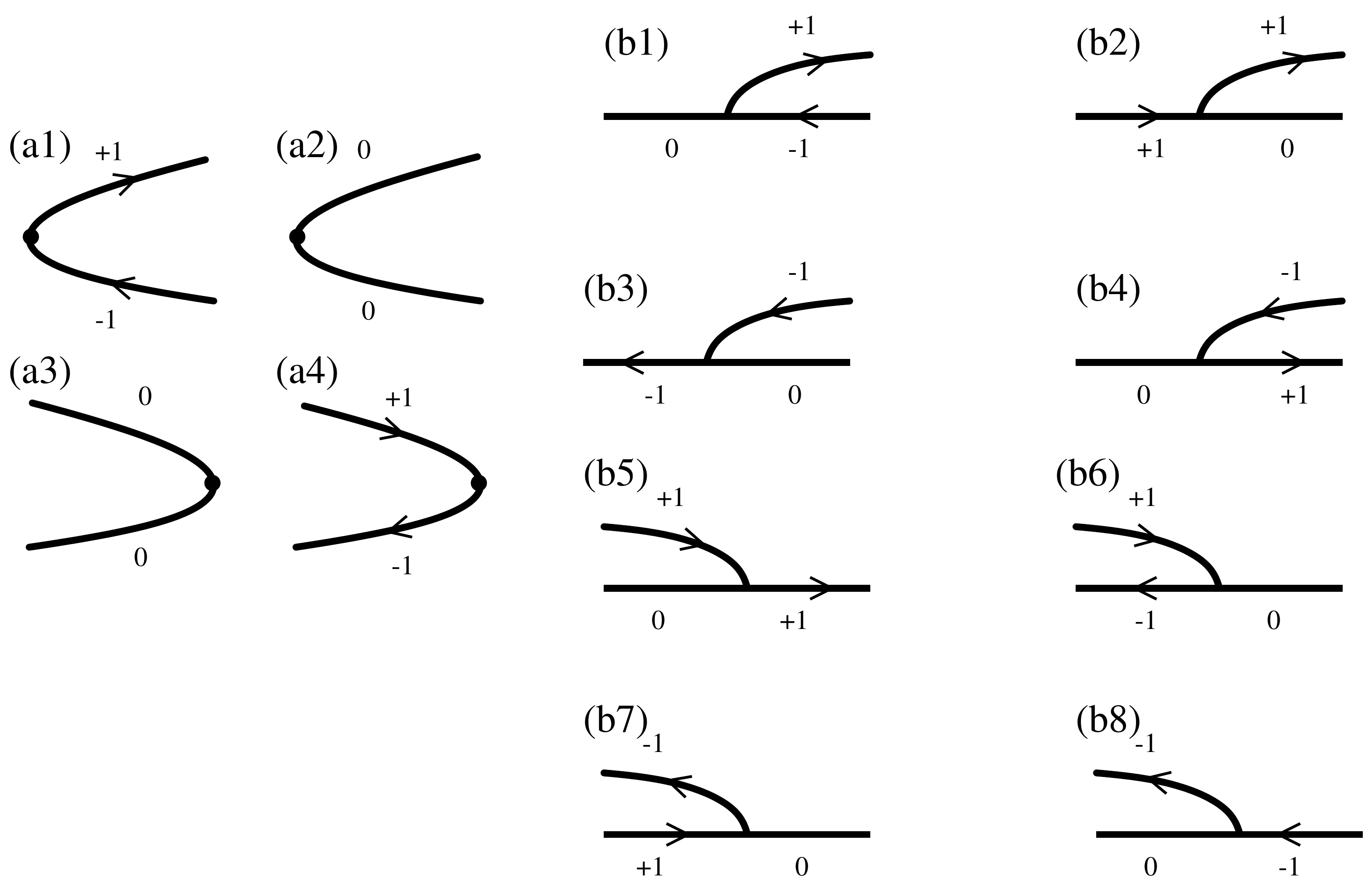}
\caption{ \label{fig:sadd}Generic Bifurcations: This figure depicts all possible generic (a)
 saddle-node and (b) period-doubling (or period-halving) bifurcations
 with their orbit indices, the numbers near each segment. Each point denotes an orbit. In this symbolic representation, the horizontal axis is the parameter, and the vertical axis is $PO(F)$. The arrows have the same meaning as in Figure~\ref{fig:snake} but segments of flip orbits are indicated here by having no arrows (and by having orbit index $0$.
 The one-dimensional quadratic map $\lambda - x^2$ has only bifurcations of types (a1) and (b2).}
\end{center}
\end{figure}

\section{Theorem of cascades from boundaries\label{sec:abstracttheorem}}

\subsection{Oriented arcs entering or exiting regions\label{sec:boundedregions}}

We now describe the restriction of oriented arcs to a region $U$ with
a bounded parameter range.

\begin{hypothesis}[Orbits near the boundary]\label{hyp:3}
  Let $F$ satisfy Hypothesis~\ref{hyp:2}. Let $\lambda_0<\lambda_1$, and let
  $U= [\lambda_0,\lambda_1] \times \mathfrak{M} $ and 
$ \partial U=\{\lambda_0,\lambda_1\} \times \mathfrak{M}$.  
Assume that all orbits in $ \partial U$ are hyperbolic. Assume that all
  orbits in $U$ are contained in a compact subset of $\mathfrak{M}$.
\end{hypothesis}

\begin{definition}
Assume Hypothesis~\ref{hyp:3} and its notation. 
Let $p \in PO_{nonflip}(F)$ be a hyperbolic
  orbit in $ \partial U $. If $p$ is oriented ``into'' the region
  $U$ by index orientation, then it is called an {\bf entry orbit} of
  $U$. That is, $p$ is an entry orbit if either $\lambda = \lambda_0$ and
  $\phi(p)=+1$, or $\lambda= \lambda_1$ and $\phi(p)=-1$. Otherwise, 
  it is called an {\bf exit orbit} of $U$.

A cascade is said to be {\bf essentially in $U$} if all but a finite number of its bifurcation orbits are in $U$.

\end{definition}

\begin{theorem}[Cascades from boundaries]\label{thm:main}
 Assume Hypothesis~\ref{hyp:3}. 
Let $IN$ be the
 set of entry orbits. Let $OUT$ be the set of exit orbits. 
Assume that $IN$ contains
 $K$ elements, and $OUT$ contains $J$ elements. We allow one but not
 both of the sets to have an infinite number of elements.

\begin{enumerate}
\item[($C_K$)] 
If $K<J$,
then all but possibly $K$ orbits in $OUT$ are contained in distinct components, each of which contains a cascade that is essentially in $U$.

Likewise, if $J<K$, 
then all but possibly $J$ orbits in $IN$ are contained in distinct components, each of which contains a cascade that is essentially in $U$.

\item[($C_0$)]
 If $J = 0$ or $K = 0$, then 
the nonflip orbits of $\partial U$
are in one-to-one correspondence with the components that intersect the boundary. Each of these components has one cascade that is essentially in $U$. \end{enumerate}
\end{theorem}

\begin{proof}
 Assume Hypothesis~\ref{hyp:3}. For simplicity, we specify that
 $J<K$. The proof of the other case is similar. 

Let $q$ be an orbit of $IN$. Then $q$ is an entry orbit and it lies in a 
 component $Q$, which by
 Theorem~\ref{thm:OneManifold} is a one-manifold. 
Let $\sigma: X \to Q$ be the index orientation for $Q$.
Starting at $q$ and following the component in the forward direction (using the index orientation), 
$\sigma(s)$ initially enters the interior of U. 
If it leaves $U$, it does so through an exit orbit.
Let $exit(q)$ denote the first such exit orbit encountered.  
Distinct entry orbits $q$ that exit yield distinct $exit(q)$ of which there are at most $J$. Hence all but $J$ entry orbits are in components that do not leave $U$ for increasing $s$. Such a component must be an open arc and the component must have a bounded end. By Proposition~\ref{prop:bddArcs}, it contains a cascade, a cascade that is essentially in $U$. Hence ($C_K$) is proved. 

If $J = 0$, then each entry orbit is in a component that crosses the boundary only once. Furthermore any component that crosses the boundary must do so at an entry orbit. Each such component has an end in $U$ and so has a component that is essentially in $U$, proving ($C_0$). 
\end{proof}

\subsection{Related Results}\label{sec:related}

Our abstract results in the previous section build on the work in~\cite{yorke:alligood:83} and Franks~\cite{franks:85}. 
We now compare our results to these two previous results.

The papers~\cite{yorke:alligood:83, yorke:alligood:85} proved the
existence of cascades of attracting periodic points for area
contracting maps and for elliptic periodic points for area preserving
maps, in a particular case of Theorem~\ref{thm:main} Part $C_0$, without
assuming genericity. Whereas our current result applies to
parametrized maps with a large number of unstable dimensions, their
result only considered maps with at most one unstable dimension. The
existence of attractors relies on having at most one unstable
dimension, because this implies that there are no Hopf
bifurcations. For parametrized maps with more than one unstable
dimension, cascades do not in general contain attractors. Both results
involve snakes in the generic case, followed by smooth convergence
arguments to show the general case. These convergence arguments no
longer apply when there is more than one unstable dimension.

The {\bf Morse index} is the number of unstable eigenvalues. If the Morse index is even, the orbit index is either 0 or $+1$. If the Morse index is odd, the orbit index is either 0 or $-1$. 

Franks~\cite{franks:85} proves there are cascades under the following conditions. 
Let d be an odd integer.
Assume that for every non negative integer k, every orbit of period $2^kd$ has a Morse index 
with the same ``parity'' (all are odd or all even) at $F(\lambda_0,\cdot)$, and the opposite parity at $F(\lambda_1,\cdot)$. In
our notation, this corresponds approximately to saying that on the boundary $F(\lambda_0,\cdot)\cup F(\lambda_1,\cdot)$, all orbits are entry orbits (or alternatively all are 
exit orbits); we ignore flip orbits. Our theorem relaxes this condition; we only 
assume that the numbers of entry and exit orbits differ. 

Franks' proof uses the Lefschetz trace formula, which allows the
smoothness of $F$ to be relaxed. $F$ is assumed to be a continuous
parametrized map, where $F(\lambda_0,\cdot)$ and $F(\lambda_1,\cdot)$ are smooth maps. However,
this lack of smoothness has implications. The theorem does not give
information about the bifurcations, only assuring that the component
of $PO(F)$ containing the original hyperbolic orbit of period
$2^rd$ ($d$ odd) contains flip orbits of period $2^kd$ for all
$k \in N$ on the boundary $F(\lambda_0,\cdot) \cup F(\lambda_1,\cdot)$. The theorem says nothing
about the way these orbits bifurcate. A cascade is usually
viewed as a sequence of events with some separation, but in
the context of Franks, a portion of the cascade can occur at a single parameter value. For
example, in non-generic maps with dimension larger than one, two
eigenvalues can simultaneously bifurcate through $-1$. Thus a fixed
point can bifurcate to a period-four orbit, missing a bifurcation through period two. A 
more extreme case of this phenomenon is shown in the following example using a slight adaptation of Franks' result.
The example is a one-dimensional but non-smooth map, in which an entire generalized cascade occurs at one parameter value. 

\begin{example}[The parametrized tent map] {\em Consider the
    parametrized map consisting of tent maps of slopes of absolute
    value $\lambda$, with $\lambda$ increasing to the standard value
    of $2$.  Orbits of all the periods $2^k$ (for $k \ge 0$) appear at
    exactly $\lambda=1$, or more precisely, each exists for all
    $\lambda > 1$. In fact, for every $k>2$, the period-$k$ cascades
    appear in this manner. All of these are cascades in the sense of
    Franks' theorem.  }
\end{example}

\section{Cascades for new classes of functions }\label{sec:examples}

In this section, we describe a number of examples of classes of parametrized maps, each with an infinite
number of cascades. 

\subsection{ Cascades for parametrized maps with horseshoes }\label{applications}

In this section we show under certain restrictions that the creation of a Smale horseshoe in dimension two implies 
the existence of infinitely-many cascades. This is a significant
generalization of the work of Yorke and
Alligood~\cite{yorke:alligood:83}, who showed that if as a parameter
changes, a map develops a Smale horseshoe in a very specific manner,
then there are cascades.
Here we make no assumptions about how the
horseshoe is created and instead make only mild assumptions about
other orbits.

Rather than defining a Smale horseshoe, we prefer to simply
state the properties we use, namely $(S_1)$--$(S_2)$ in Corollary~\ref{cor:SmaleHorseshoe} below. Our assumption
$(S_1)$ that $F_1$ has no bifurcation orbits seems reasonable
since generically there are at most countably many bifurcation orbits
and hence, at most countably many values of $\lambda$ at which there
are bifurcation orbits.

\begin{corollary}[Creating a Smale horseshoe]\label{cor:SmaleHorseshoe}
Let the dimension of ${\mathfrak M}$ be $2$.
Assume Hypothesis~\ref{hyp:2}.
Let $W = [\lambda_0, \lambda_1] \times {\mathfrak M}$.
Assume the following:
\begin{enumerate}
\item[($S_0$)] $F_0 = F(\lambda_0, \cdot)$ has at most a finite number of saddle orbits.
\item[($S_1$)] $F_1 = F(\lambda_1, \cdot)$ has at most a finite number of orbits 
that are either attractors or repellers, and all its orbits are hyperbolic. 
\item[($S_2$)] $F_1$ has infinitely many nonflip saddle orbits.
\item[($S_3$)] There is a compact subset of $W$ that contains all the orbits in $W$. 
\end{enumerate}
Then there are infinitely many cascades whose components have a bounded end in $W$. 
\end{corollary}

Concerning ($S_1$), it is conjectured that it is always possible to choose that $\lambda_1$ so that there are only finitely many attracting or repelling orbits at
that value. 
Even though diffeomorphisms with infinitely many coexisting sinks are Baire generic (as Newhouse showed \cite{Newhouse:74}), it is conjectured that they have ``probability zero'' in the sense of prevalence. See Gorodetski and Kaloshin \cite{Gorodetski-Kaloshin:2007} for recent partial results in this direction. The first results in this direction were much earlier in \cite{Ted-Yorke:86} and \cite{Nusse-Ted:92}. 

We may thus plausibly expect a typical parametrized map to have finitely many attractors for almost every parameter value.
If this property is true, we can apply it to inverses of maps to conclude that there are only finitely many repellers for almost every parameter value. Hence we can plausibly assume that for almost every parameter value there are finitely many attractors and repellers. Even if the conjecture is false, our assumption is true for many systems. 

\bigskip

\begin{proof} The saddles in $S_2$ are all entry orbits and they are infinite in number. There are at most a finite number of exit orbits. Hence the result follows from Theorem~\ref{thm:main} part $C_K$.
\end{proof}

Dynamical systems that satisfy conditions $(S_0), (S_2)$, and $(S_3)$ for some $\lambda_0$ and $\lambda_1$ are plentiful. These give evidence of satisfying $(S_1)$, but a rigorous check of this condition is generally difficult at best and not practical. The following processes are examples. We have chosen specific parameter values, though the phenomena described are seen over a wide range of parameter intervals.

{\bf The Ikeda map} models the field of a laser cavity \cite{Hammel-Jones-Moloney}. 
For $z \in {\mathfrak M} = \mathbb C$, the complex plane, the map is
\[F(\lambda, z) = \lambda + 0.9 \; z \; e^{i \{ 0.4 - 6.0/(1 + |z^2|) \}} \]
At $\lambda = 0$ there appears to be a globally attracting fixed point. 
At $\lambda = 1.0$ we observe numerically a global chaotic attractor with a positive Lyapunov exponent and homoclinic points and one attracting fixed point and no repellers.

\bigskip

{\bf The Pulsed Rotor map} with $(x,y) \in {\mathfrak M} = S^1 \times R$ is
\[F(x,y) = ((x + y) (\mod 2\pi ), 0.5 y + \lambda \sin(x+y)). \]
For $\lambda = 0$, there is a saddle fixed point and an attracting fixed point that attracts everything except for the stable manifold of the saddle. For $\lambda = 10$ we observe a chaotic attractor and a fixed point with a transverse homoclinic point.

\bigskip

\noindent
{\bf Geometric Off-On-Off Chaos.} 
In his investigations of the Lorenz
system, John Guckenheimer \cite{Guckenheimer:1976} realized that we were not close to a
rigorous understanding of the system, so he introduced a {\it
 geometric} Lorenz model. He gave this model well defined geometric
properties that the Lorenz system appeared to have. Those properties had not been rigorously established at that time. 
In that way he could see what these
conjectured properties implied. We follow a similar path with the
double-well Duffing equation and forced damped pendulum. See Figure~\ref{fig:duffing}. Both are 
periodically forced with period $2\pi$. Therefore their time-$2\pi$ maps,
denoted by $F(\lambda,\cdot)$ or $F(\lambda, u, du/dt)$ are diffeomorphisms on
$R^2$. 

Our numerical investigations
strongly suggest that for the double-well Duffing equation 
(1) there are a number of intervals in the
parameter range where F has a globally attracting fixed point. These
intervals are centered near the values $\lambda \in \{1.8, 20., 73.,
 175., 350.\}$. 

Furthermore (2) between any two consecutive values in
this set, there is a $\lambda$ for which (some iterate of) the time - $2\pi$ map appears
to have a Smale horseshoe. Actually it seems to have a chaotic
attractor. Zakrzhevsky~\cite{zakrzhevsky:2008} provides many insights into the dynamics of a double-well Duffing equation though his version uses a restoring force of $u^3 + u$ instead of our choice of $u^3 - u$.

We formalize these properties in the following definition.

\begin{definition}
Assume the dimension of ${\mathfrak M}$ is 2.
 We call a map $G : R \times {\mathfrak M} \to {\mathfrak M}$ a {\bf Geometric
 Off-On-Off-Chaos} map 
if it satisfies the following
 properties.
\begin{enumerate}
\item[($D_0$)] There are values $\Lambda_1 < \Lambda_3$ such that $G_1 =
 G(\Lambda_1, \cdot)$ and $G_3 = G(\Lambda_3, \cdot)$ each have
at most a finite number of orbits, whose total is $k$.
\item[($D_1$)] There is a $\Lambda_2 \in (\Lambda_1, \Lambda_3)$ for
 which $G_2 = G(\Lambda_2, \cdot)$ has at most a finite number of
 orbits that are attractors or repellers, and all of its orbits
 are hyperbolic.
\item[($D_2$)] $G_2$ has infinitely many nonflip saddle orbits.
\item[($D_3$)] There is a compact subset of $W= [\Lambda_1, \Lambda_3] \times {\mathfrak M}$ that contains all of the orbits in $W$. 
\item[($D_4$)] $G$ satisfies Hypothesis~\ref{hyp:2}.
\end{enumerate}
\end{definition}

\begin{corollary}\label{cor:OffOnOffChaos}
 Assume $G$ is a geometric on-off-on-chaos map. Then there are infinitely many pairs of
 cascades in $W= [\Lambda_1, \Lambda_3] \times {\mathfrak M}$, where the two cascades
 of each pair are in the same component of
 $PO_{nonflip}(G)$. Also there are at most $k$ unbounded cascades that have $\Lambda$ values only in $[\Lambda_1, \Lambda_3]$.
\end{corollary}

\begin{proof} Apply Corollary~\ref{cor:SmaleHorseshoe} twice 
where $W =[\lambda_0, \lambda_1]$, 
once to $W_1 =
 [\Lambda_1, \Lambda_2] \times {\mathfrak M}$ and once to $W_2 = [\Lambda_2, \Lambda_3] \times
 {\mathfrak M}$. It follows that there are infinitely many bounded components
 of $PO_{nonflip}(G)$ that have one bounded end in $W_1$ and the
 other in $W_2$. Any such component contains a pair of cascades, proving the first assertion.

If $C$ is an unbounded cascade that has $\Lambda$ values in $[\Lambda_1, \Lambda_3]$, then its component is unbounded and must have an orbit whose $\Lambda$ coordinate is either $\Lambda_1$ or $\Lambda_3$. Since the map has only $k$ such orbits, there can be at most $k$ such unbounded cascades.

\end{proof}

{\bf The forced damped pendulum.}
The time-$2\pi$ map for 
\[\frac{d^2\theta}{dt^2} + 0.3 \; \frac{d\theta}{dt} + \sin \theta = \lambda \cos t \]
strongly appears to yield a geometric off-on-off-chaos map.  We investigate the time-$2\pi$
map on ${\mathfrak M} = S^1 \times R$; that is, the first variable is
$\theta~(\mod~2\pi)$ and the second is $d\theta/dt \in R$.  There is a
symmetry about $\lambda = 0$: For parameters $\lambda$ and $-\lambda$
the system has the same dynamics.  At $\lambda = 0$, there are only
two periodic orbits, both fixed points, an attractor and a saddle, and
we observe numerically a global chaotic
attractor with a positive Lyapunov exponent and homoclinic points at $\lambda = 2.5$. 
For
$\lambda \ge 10$, as at $\lambda = 0$, the two fixed points are the only
orbits.  Due to the friction term $0.3 \; d\theta /dt$, the
orbits must lie in a compact subset of ${\mathfrak M}$ for
$\lambda \in [-10,10]$.  This map then appears to be a geometric off-on-off-chaos map with
either $\Lambda_1 = 0$, $\Lambda_2=2.5$, and $\Lambda_3 = 10$, or by
symmetry, with $\Lambda_3=0$, $\Lambda_2=-2.5$, and
$\Lambda_1=-10$. Assuming these numerical observations are valid, each
cascade must have its $\lambda$ values lie entirely in either $(0,10)$
or $(-10,0)$. There are at most $k=4$ unbounded cascades and an
infinite number of bounded pairs of cascades.

\subsection{Large-scale perturbations of a quadratic map}\label{subsec:e2}

In this subsection, we state Corollary~\ref{cor:1D-quadratic} for maps
of the form $\lambda-x^2 + g(\lambda,x)$. We first give a method for
counting cascades of period $k$ for each $k$.

{\bf The number $\Gamma(1,k)$.} We will describe the number of cascades in terms of the tent map,
\begin{equation}
T(x) = \begin{cases}
2x & \mbox{ for } x \in [0, 0.5]\\
2(1-x) &\mbox{ for } x \in (0.5, 1].
\end{cases}
\end{equation}

An equivalent formulation of {\bf nonflip} for this map is the following. An orbit is nonflip if it has an even number of points in $(0.5, 1]$. If the number of points is odd, it is a flip orbit.
Nonflip period-k orbits are the orbits whose derivative satisfies $\frac{d}{dx} (T^k)(x) > +1$. The derivative is $< -1$ for flip orbits and no orbits have derivative in $[-1, 1]$.
Define $\Gamma(1,k)$ to be the number of period-$k$ nonflip orbits of the tent map $T$. (The entry ``1'' refers to the dimension of $x$.) We give a general formula for $\Gamma$ in Section~\ref{subsec:e2}.

We are now ready to state a result:

\begin{corollary}[A large-scale perturbation of the parametrized quadratic map]\label{cor:1D-quadratic}
Assume that $F:R \times R \to R$ has the form
\begin{equation} F(\lambda,x)=\lambda-x^2+g(\lambda,x),\end{equation}
where $g:R \times R \to R$ is $C^\infty$.
Assume that there is a $\beta>0$ such that for all $\lambda$ and $x$, $|g(\lambda,0)|<\beta$, and $|\frac{\partial g}{\partial x} (\lambda,x)| < \beta$. 
 Then for a residual set of $g$,
for each positive integer $k$, the number of unbounded period-$k$ cascades is
$\Gamma(1,k)$,
which is the same as for $\lambda-x^2$.
\end{corollary}

{\bf Unbounded cascades cannot be destroyed.} Here is an interpretation of the above corollary. Let $\gamma \ge 0$ be sufficiently large that all orbit bifurcations and all attractors of $\lambda-x^2$ lie in the arbitrarily large square
\[S = \left\{ (\lambda, x) \in [-\gamma, \gamma]\times [-\gamma, \gamma] \right\}. \] We can choose $g$ from the residual set in the corollary so that everywhere on $S$, $F(\lambda, x) = 0$.
Hence all the cascades and all of the bifurcations have been eliminated from the arbitrarily large bounded set $S$. The corollary simply says then that the cascades still exist but they must exist outside of $S$. Such $g$ can be chosen so that $F$ is generic.

In addition to unbounded cascades, there may be bounded pairs of cascades, as is seen in Figure~\ref{fig:boundedcascades}.

\bigskip

One key fact for the proof is that for
$\lambda=\lambda_H$ sufficiently large, $F(\lambda_H,\cdot)$ is a
two-shift horseshoe map, as defined in the next definition.

\begin{definition}[ Two-shift horseshoe map in dimension 
 one]\label{d:horseshoe}
 We refer to a $C^1$ one-dimensional function $G: R \to R$ as a {\bf two-shift horseshoe map} when it has the following properties:
\begin{enumerate}

\item there is a closed interval
 $J$ and two non-empty disjoint intervals $J_1 \subset J$
 and $J_2 \subset J$ such that  $G(J_1) = G(J_2) = J$;

\item $G(x) \in J$ implies $x \in J_1 \cup J_2$; 

\item $G'(x) < -1$ for $x \in J_1$ and $G'(x) > 1$
 for $x \in J_2$.
\end{enumerate}
\end{definition}

Assume that $F$ is a large-scale perturbation of a parametrized quadratic map as
in Corollary~\ref{cor:1D-quadratic}, and $\lambda_H$ is sufficiently large that
$F(\lambda,\cdot)$ is a two-shift horseshoe map
for $\lambda \ge \lambda_H$.
Let $y$ be a nonflip
orbit for $F(\lambda_H,\cdot)$, and let $C$ be the component
containing $y$. We show in
Proposition~\ref{prop:minPerQuadraticMap} below that the minimum period of the orbits
contained in the component $C$ is equal to the period of $y$, which is the stem period of the component. In
the proof of Corollary~\ref{cor:1D-quadratic} below, we show that $C$ contains
an unbounded cascade, and no other nonflip orbit at
$\lambda=\lambda_H$ is contained in $C$. That is, there is a
one-to-one correspondence between the unbounded cascades of $F$ and
the orbits of nonzero index in the two-shift horseshoe. Combining these two
results, we conclude that there is a one-to-one correspondence between
the unbounded {\it period-$M$ cascades} and the {\it period-$M$
 orbits} of nonzero index in the two-shift horseshoe. The following
definition describes this relationship more precisely.

\begin{definition}[Stem period]
Assume Hypothesis~\ref{hyp:2}. Let $C$ be a cascade in an unbounded component $Q$.
Assume there is a compact subset $B$ of $ R \times {\mathfrak M}$ such that all orbits of $C$ that do not lie in $B$ have the same period.
We call that period the {\bf stem period} of $C$. 
\end{definition}

\begin{remark}
{\em {\bf The period of an unbounded cascade versus its stem period.}
Proposition~\ref{prop:minPerQuadraticMap} below shows that for large-scale perturbations of the quadratic parametrized map,
each nonflip period-$k$ orbit of the two-shift horseshoe map $F(\lambda_H,\cdot)$ is contained in a component with least period $k$. That is, in the large-scale perturbed quadratic case, the stem period is equal to the period of the cascade. The proof relies on the parametrized map having only one spatial dimension, in addition to the fact that for sufficiently large $\lambda$, the map is a typical two-shift horseshoe map.

It is always the case that if a period-$k$ orbit is contained in a cascade, then the cascade is a period-$M$ cascade, where $k/M$ is a power of two. Thus for odd $k$, a period-$k$ cascade can be of no lower period than its stem period, since there is no period-halving
bifurcation for an odd period. Therefore, 
if the stem period of a cascade is odd, it is equal to the period of the cascade. 
 We conjecture that in general for even $k$, the stem period need not be equal to the period of the cascade, such as for the classes of non-quadratic parametrized maps in Corollary~\ref{cor:cubic} or the classes of higher-dimensional parametrized maps in Corollary~\ref{cor:QuadSys}.
 }
\end{remark}

We now proceed with the proposition. This is the only part of the proof of Corollary~\ref{cor:1D-quadratic} which uses the one-dimensionality of the phase space.

\begin{proposition}[Stem Period = Minimum Period]\label{prop:minPerQuadraticMap}
Assume that $F: R \times R \to R$ satisfies Hypothesis~\ref{hyp:3} (with dimension $N = 1$), and that for all sufficiently large $\lambda$, $F(\lambda_H,\cdot)$ is a two-shift horseshoe map.
Let $G = F(\lambda_H, \cdot)$, where $\lambda_H$ is chosen such that $G$ is a two-shift horseshoe map.
Let $y$ be a period-$M$ orbit for $G$ with nonzero orbit index. Let $C$ be the component containing $(\lambda_H, y)$. Then $M$ is the minimum period of an orbit in $C$.
\end{proposition}

\begin{proof}
The argument proceeds by showing that if the minimum period were less
than $M$, then $y$ would be a flip orbit, contradicting our assumption
that it has nonzero index.

 Assume the contrary, and that the component $C$ of $y$ 
contains an orbit of period less than $M$. Then $M$
 is even, and $C$ contains an orbit of period $M/2$. Since $C$ is a
 one-manifold, there is a continuous $\sigma : [0,1] \to
 PO_{nonflip}(F)$ such that (a) $\sigma(0) = (\lambda_H, y)$,
 (b) $\sigma(1)$ is an orbit of period $M/2$, and (c) the period of
 $\sigma(s)=(\lambda(s),y(s))$ is at least $M$ for $s \in [0, 1)$.
 As $s$ increases, the set $y(s) \subset R$ can undergo period doublings or
 halvings. However, no other changes can occur, since we have assumed generic
 bifurcations.

 Note that since y is an orbit in dimension one, there is a well-defined
 ordering of the points in the orbit $(\lambda(s),y(s))$, and
 this order is preserved under variations of $s$, except at a periodic
 orbit bifurcation. We now show that the ordering of points in an orbit is preserved at
such bifurcations as well. We first formalize the ordering of points in an orbit.

{\bf Splitting $\sigma(s)$ into $M$ subsets ${S_k(s)}$.} Define the
subsets $S_1(s), \dots, S_M(s)$ of C for $s \in [0,1)$ as follows,
with each 
$S_k(s) \subset \{\lambda(s)\} \times R$.
For $s=0$, each contains one point in $R \times R$. Let $S_1(0)$ consist
of the point $(\lambda_H, y_1)$ where $y_1$ is the smallest (i.e., most negative) point of the orbit $y$,
$S_2(0)$ of the point $(\lambda_H, y_2)$, where $y_2$ is the next smallest point,
etc., up to $S_M(0)$, where $y_M$ is the largest point. For each $k$, let
$S_k(s)$ be the continuation of $S_k(0)$ in $C$. For each $s \in
[0,1)$, if $\sigma(s)$ has period $q$, then $q/M$ is a power of two, and
each $S_k(s)$ contains $q/M$ points. Note that for
each $k$, $S_k$ varies continuously with respect to $s$ in the Hausdorff metric.

{\bf Permutations and natural pairs.} Let $\pi$ be the permutation
on $\left\{ 1,\dots,M \right\}$ defined by $F(S_k(0)) =
S_{\pi(k)}(0)$. Then $F(S_k(s)) = S_{\pi(k)}(s)$ for all $s \in [0,
1)$. We will call a pair of integers {\bf natural} if they are
consecutive and the larger one is even (i.e. $\left\{ 1,2 \right\},
\left\{ 3,4 \right\}, \dots$). Hence $\pi$ maps natural pairs to
natural pairs. 

There is a strong relationship between period-halving bifurcations and
natural pairs: Assume there is a period-halving bifurcation of
$\sigma(s)$ at $s_*<1$. By the above discussion, prior to the
bifurcation, we write $\sigma(s)=(\lambda(s),\left\{ y_1(s), \dots,
 y_q(s) \right\})$, where we always order the orbit such that
$y_1(s) < y_2(s) < \dots < y_q(s)$. Note that $q$ is even, since it is about to
go through a period-halving bifurcation. By the constraints of
the ordering on the real line, at $s_*$, there is a merging of $y_i(s)$
and $y_{i+1}(s)$ for every natural pair $\left\{ i, i+1 \right\}$ of
integers in $\left\{ 1, \dots, q \right\}$. Further, if the orbit
period doubles and then halves, the pairs created in the
doubling are the same as the pairs that merge back together when
halving.

{\bf $S_k(s)$ are disjoint for $s<1$.} The sets $S_k(s)$ do not merge for $s < 1$. If such a merging occurred,
it would have to occur at a period-halving bifurcation of the orbit with
the period going from $q$ to $q/2$, where $q/2 \ge M$. However,
the largest point in $S_1(s)$ is an even numbered point, which would
merge with the smallest point in $S_2(s)$, an odd numbered point. This
would not be a natural pair, which is a contradiction to the argument
above. Thus the $S_k(s)$ do not merge for $s < 1$.

At $s=1$, the period drops from $M$ to $M/2$. The $M$ sets $S_k(s)$
contain a single point as $s$ approaches 1, and at $s=1$ the sets
merge in natural pairs. That is, $S_i$ and $S_{i+1}$ merge exactly
when $\left\{ i,i+1 \right\}$ is a natural pair. To preserve
continuity, this implies that for every $0 \le s \le 1$,
$F(\lambda(s),\left\{ S_i,S_{i+1} \right\} ) = \left\{ S_j,S_{j+1}
\right\}$, where both $\left\{ i, i+1 \right\}$ and $\left\{ j, j+1
\right\}$ are natural pairs. That is, $F$ maps natural pairs of sets to
natural pairs of sets.

Let $J_1$ and $J_2$ be the expanding intervals for $G$, which are
guaranteed to exist since $G$ is a typical two-shift horseshoe map. If
$y_{2k-1}(0)$ and $y_{2k}(0)$ are either both in $J_1$ or both in $J_2$,
then, $|y_{2k-1}(0)- y_{2k}(0)| < |G(y_{2k-1}(0))- G(y_{2k}(0))|$.
Since $G$ is a permutation of the natural pairs in $[y(0)]$, the above
inequality cannot hold for all $k$, so there must be one pair which
has a point in $J_1$ and another in $J_2$. This pair is the only such
pair, since points are linearly ordered. Hence there is an odd
number of points of the $[y]$ orbit in $J_1$, implying that there is
an odd number of points of the orbit at which the derivative $\frac{dG}{dx}$
is negative. Therefore $\frac{d}{dx}(G^{M})(y_i) < -1$ at each point of the
orbit. This implies the orbit index of $(\lambda_H,y)$ is zero, and
$y$ is a flip orbit. This contradicts our assumption that $y$ has a nonzero orbit index.
Thus the least period in $C$ is $M$.
\end{proof}

We now proceed with the proof of Corollary~\ref{cor:1D-quadratic}.

\begin{proof} (of Corollary~\ref{cor:1D-quadratic})
Assume the hypotheses of the corollary.
The proof of this corollary follows from showing that the following hold:
\begin{enumerate}
\item For all sufficiently small $\lambda_0$, there are no orbits.
\item For sufficiently large $\lambda_1$, the maximal bounded invariant set of $F(\lambda_1, \cdot)$ is
conjugate to a shift map on $2$ symbols.
To show that, we prove that on this set, the map is expansive. (We say that a map $F: R \times R \to R$ is {\bf expansive} on a set if at every point of the set $|F_x| > 1$.)

\end{enumerate}

After these facts are proved, Theorem~\ref{thm:main} part $C_0$ will complete the proof.

We now show that these facts hold.
Let $s = \sqrt{\lambda_1}$ (where $\lambda_1 > 0$) and $J = [-2s,+2s]$. Let $A_1=[-2s,-s/2]$, and $A_2=[s/2,2s]$, and let $A=A_1 \cup A_2$.
We show the following estimates for $\lambda_0$ chosen sufficiently small, and $\lambda_1$ chosen sufficiently large:
\begin{enumerate}
\item[(a)] At $\lambda_0$, there are no orbits, and orbits of $F$ are strictly monotonically decreasing.
\item[(b)] At $\lambda_1$, $|\partial F/\partial x|>1$ for any point in $A$.
\item[(c)] At $\lambda_1$, the set $F(\lambda_1,R \setminus A)$ contains no points of $J$.
\item[(d)] There is an interval $B$ in the interior of $J$ such that for all $\lambda<\lambda_1$, each orbit is contained in $B$.
\item[(e)] At $\lambda_1$, $F(\lambda_1, A_i)$ maps diffeomorphically across $J$, where $i=1$ or $2$.
\end{enumerate}

The proof of 1 and 2 follows from these five estimates as follows: It follows from (b)-(e) that at $\lambda_1$, $F$ on the maximal invariant
set in $J$ is topologically conjugate to the full two shift on symbols
associated with $A_1$ and $A_2$. Specifically, all orbits are
hyperbolic. It follows from (d) that there are no orbits on or near
the boundary of $J$ for all $\lambda$ between $\lambda_0$ and
$\lambda_1$. It follows from (a) that there are no orbits on or
near $\lambda=\lambda_0$. Let $C$ be an arc in $PO_{nonflip}(F)$
in region
$U= [\lambda_0,\lambda_1] \times \mbox{int} (J)$ 
with index orientation $\sigma: I \to C$
(where $I$ is the interval $[0,1]$). If either of $A=\sigma(0)$ or $Z=\sigma (1)$ is on
the boundary of $U$, then it must be on the piece of the boundary
$\lambda=\lambda_1$, since there are no orbits near any other portion of the boundary.
All orbits on $\lambda=\lambda_1$ are hyperbolic orbits, implying that all orbits near the boundary have the same orbit index as the limit orbit on the boundary. Furthermore, they are one-dimensional unstable orbits, implying that they are either flip orbits or that the orbit index is equal to
$-1$. This implies that there are no entry orbits for $U$ ($IN$ is empty). 
Each nonflip orbit of $F(\lambda_1,\cdot)$ is an 
entry orbit. Thus by 
Theorem~\ref{thm:main} part $C_0$,
$F$ has exactly one cascade containing each orbit with nonzero orbit index at $\lambda=\lambda_0$.

We now verify the five estimates listed above.

\begin{itemize}

\item[Proof of (a)] Let $\lambda_0 < -(1+6 \beta + \beta^2)/4,$, where $\beta$ is the bound given in Corollary~\ref{cor:1D-quadratic}. This
 condition guarantees that for all $x$, $\lambda_0 - x^2 + \beta (1+|x|) -x<0$. Therefore
\[ F(\lambda_0,x) -x < \lambda_0 - x^2 + \beta (1+|x|) -x <0. \]
Thus $F$ lies below the diagonal, so it is strictly monotonically decreasing. Therefore, there are no periodic points when $\lambda=\lambda_0$.

\item[Proof of (b)] Assume that $\lambda_1>4 (\beta+2)^2$.
 We now show that $|\partial F/\partial x| >1$ as long as $s/2 \le |x| \le 2s$:
\[s^2 = \lambda > 4 (\beta+2)^2> 4 (\beta+1)^2,\] implying
that \[2|x| \ge |x| \ge s/2 > \beta+1 > |\partial g/\partial x| + 1.\]
Therefore \[ \left| \frac{\partial F}{\partial x} \right| = \left| -2x
 + \frac{\partial g}{\partial x} \right| \ge |-2x|- \left|
 \frac{\partial g}{\partial x} \right| > 1.\]

\item[Proof of (c)] We now show that if $|x|\le s/2$, then $F(\lambda_1,x)>2s$.
Assume that $|x| \le s/2$. The conditions on $g$ imply that for all $\lambda$
and $x$, $|g(\lambda,x)| < \beta + \beta |x|$. Thus
\begin{eqnarray*}F(\lambda_1,x)-2s &>& \lambda_1 -x^2 -\beta (1+|x|) -2s\\
&\ge& s^2 - \frac{s^2}{4} - \beta \left( 1+\frac{s}{2} \right) -2s
= s \left( \frac{3s}{4} - \frac{\beta}{2} - 2 \right) - \beta \\
&\ge& s \left( \frac{3 (\beta+2)}{2} - \frac{\beta}{2}-2 \right) - \beta \\
&\ge& 2 (\beta+2) (\beta+1) -\beta >0.
\end{eqnarray*}

We now show that if
$|x|>2s$, then $F(\lambda_1,x)<-2s$.

Assume that $|x| \ge 2s$. Then
\begin{eqnarray*}
F(\lambda_1,x)+2s &<& \lambda_1 - |x|^2+ \beta (1+|x|)+2s \\
&\le& s^2+2s +|x| (-|x|+\beta) + \beta \\
&\le& s^2+2s+2s (-2 s + \beta) +\beta \\
&\le& -3 s^2+2s ( \beta+1) + \beta \\
&=& s(-3 s + 2(\beta+1))+\beta \\
&\le& s (-6 (\beta+2) +2(\beta+1))+ \beta \le 2(\beta+2) (-4 \beta-10) +\beta <0.
\end{eqnarray*}

\item[Proof of (d)]
Let $\{ x_1, \dots,x_k \}$ be an orbit at parameter $\lambda \le \lambda_1$. Fix $x$ to be the $x_i$ with the maximum absolute value. This implies that \[-|x| < F(\lambda,x) \le \lambda_1+\beta(1+|x|) -x^2.\] Therefore
\[0 \le \lambda_1 + \beta+|x|(\beta+1)-x^2.\] Let $\rho=(\beta+1)/2$. Then
\[0 \le (\lambda_1+\beta)+ \rho^2 -(|x|-\rho)^2.\] Thus \[|x| \le \rho+\sqrt{\lambda_1 + \beta + \rho^2}.\] This right-hand side is monotonically increasing in $\lambda_1$. Since $\rho$ and $\beta$ are fixed, for $\lambda_1$ sufficiently large, \[\rho + \sqrt{\lambda_1+\beta+\rho^2} < 2 \sqrt{\lambda_1} = 2s.\] Thus as long as $\lambda \le \lambda_1$, we have $|x|<2s$. Since $x$ is the point of the orbit with maximum absolute value, this implies that the orbit is contained in $J$.

\item[Proof of (e)] From parts (c) and (d), the image of $A_j$ lies strictly outside
 of $J$. $F$ maps $J \setminus A$ to the right of $J$ and $R \setminus J$ to the left of $J$.
By the intermediate value
 theorem, this implies that the set $F(\lambda_1,A_j) \cap J$ is
 nonempty. From (b), $|\partial F/\partial x| > 0$, so $F$ is locally a
 diffeomorphism on $A_j$. This implies that the image of the interior
 of $A_j$ is open, so the endpoints of the interval
 $C=F(\lambda_1,A_j)$ are the images of the endpoints of $A_j$. Thus
 since $F(\lambda_1,A_j) \cap J$ is nonempty, the endpoints of $C$
 are outside $J$ but not both on the same side of $J$, so $J$ is
 contained in $C$.

\end{itemize}

This completes the proof of the five estimates.
Thus we have established that all generic large-scale perturbations of quadratic parametrized maps have the same cascades. Combining this with Proposition~\ref{prop:minPerQuadraticMap} completes the proof of the result.
\end{proof}

The number of period-$k$ points for large $\lambda$ is
the same as the number of period-$k$ points for the tent map $T$ and for the full
shift on two symbols. Denote this number by $\zeta(2,k)$. We use this number to compute the $\Gamma(1,k)$, the number of nonflip orbits for the tent map.

{\bf Computing the quantity $\Gamma(1,k)$ for general $k$.}
It is straightforward to see that if $p$ is an odd prime, then there are two period-one points and $2^p$ points of period either one or $p$. Since every orbit for the tent map $T$ has a corresponding orbit which is on the left when $T$ is on the right, and on the right when $T$ is on the left, half of the period-$p$ orbits are nonflip. This implies that for an odd prime $p$, $\Gamma(1,p)=(2^p-2)/(2p)$. 
Further, for any odd integer
$k$, half of the period-$k$ orbits are nonflip orbits.
Specifically, associate a symbol sequence with an orbit in the invariant set of $T$, which we denote by
$\left\{ -1,1\right\}$, corresponding to the right and left half of
the interval respectively. The derivative is negative on the right
and positive on the left, implying that the sign of the derivative
corresponds to the representative symbol ($1$ or $-1$) in the symbol
sequence. A period-$k$ orbit of the map can be uniquely described by
its corresponding length $k$ sequence of symbols. Thus the sign of the
derivative of a period-$k$ orbit is given by the product of the $k$
symbols representing the orbit. For any sequence of the symbols $1$
and $-1$ of (least) period $k$ with product $s$, there is a
corresponding sequence with the opposite
symbols. For odd $k$, the product is $-s$. Exchanging the sign of the symbols does not change the least period, $k$.
Therefore there is a one-to-one
correspondence between flip and nonflip orbits with least period
$k$. This implies that the number of nonflip period-$k$ orbits of $F$
is $\zeta(2,k)/(2k)$ for $k$ odd, which in particular agrees with the formula for primes for any odd prime number.

For the general case, let $L(k) = \Sigma( \Gamma(1,j))$, summing over all $j < k$ for which $k/j$ is a power of $2$. Notice $L(k) = 0$ if $k$ is odd. Then
\[ \Gamma(1,k) = (\zeta(2,k)/k - L(k))/2. \]
To prove this formula, since Corollary~\ref{cor:1D-quadratic} shows that the number of nonflip orbits is the same for sufficiently large $\lambda$ for any large-scale perturbation of the quadratic parametrized maps $F$, it suffices to show that this formula gives the number of nonflip period-$k$ orbits of $\lambda - x^2$ for large $\lambda$. We show this as follows. 

Each period-$k$ saddle-node bifurcation of the 
quadratic map is connected to orbits of the following type
for each sufficiently large $\lambda$: one even orbit of period $k$ and one
flip orbit for each of the periods $2^{m}k$ for $m = 0,1,\cdots$. Hence
for each $k$ and $j < k$ for which $k/j$ is a power of $2$, each
period-$j$ cascade gives rise through a series of period doublings to
one period-$k$ orbit, and this is a flip orbit. Hence the number of
period-$k$ orbits arising in this manner is $L(k)$. The remaining
orbits arise in pairs at period-$k$ saddle-node bifurcations. Each
saddle-node gives rise to two period-$k$ segments of orbits. One is
initially an attractor but becomes a flip orbit, still of period $k$,
for large $\lambda$ and the other remains nonflip. Hence the number of
nonflip orbits is half of $\zeta(2,k)/k - L(k)$. For large $k$, $\Gamma(1,k)$
asymptotically approaches $2^k/(2k)$ in the sense that the ratio of
these two approaches $1$ as $k \to \infty$.

\subsection{Large-scale perturbations of cubic parametrized maps}

The cubic map $x^3 - \lambda x$ forms a 3-horseshoe in much the same 
way as the quadratic map form a 2-horseshoe for large $\lambda$. See Figure~\ref{fig:cubicheur}. Thus as $\lambda$ increases,
cascades occur. Nonetheless, there are some interesting features to note. For highly negative $\lambda$ there is one orbit, namely the fixed point at zero. It has orbit index 1, since the map's derivative is greater than $1$ for this orbit. As $\lambda$ increases, 0 has a period-doubling bifurcation. The period-2 orbit is symmetric about zero. The next bifurcation is a non-generic pitchfork bifurcation, namely a symmetry breaking bifurcation in which two new period-2 orbits are created. Neither orbit is symmetric about zero, (although collectively they are). Generic large-scale perturbations destroy this symmetry, resulting in a parametrized map without a pitchfork bifurcation. 

In the following corollary, we state that under any sufficiently slowly growing additive large-scale perturbation, this map has
cascades. Define the three-shift tent map $H_3:[0,1] \to [0,1]$ as the piecewise linear map with the absolute value of the slope equal to $3$
such that $H_3$ is increasing from 0 to 1 on $(0,1/3)$ and $(2/3,1)$ and decreasing from 1 to 0 on $(1/3,2/3)$. The maximal invariant set of $H_3$ is topologically conjugate to the shift map on three symbols.

\begin{corollary}[Large-scale perturbations of cubic parametrized maps]\label{cor:cubic}
 Let $F:R \times R \to R$ be of the form
 \[ F(\lambda,x)=x^3-\lambda x +g(\lambda,x).\] Assume that $g$ is a
 $C^\infty$ smooth function such that for some $\beta>0$,
 $|g(0,\lambda)|<\beta$, and for all $\lambda$ and $x$, $|\partial g
 / \partial x|< \beta |x|$. Then for a residual set of $g$, for each
 positive integer $k \ne 2^m$ for $m \ge 0$, the number of unbounded
 cascades with stem period $k$ is the same as the number of nonflip
 orbits for the map $H_3$. For all but possibly one nonflip
 period-$k$ orbit where $k=2^m$ ($m \ge 0$), there is an unbounded
 cascade component with stem period-$k$ through the orbit.
\end{corollary}

\begin{proof}
Assume that $g$ is such that $F$ satisfies genericity Hypothesis~\ref{hyp:2}.
We show the following:
\begin{enumerate}
\item[(1)] For sufficiently small $\lambda_0$, $F(\lambda_0,\cdot)$ has only one orbit, and that is an unstable fixed point.
\item[(2)] For all sufficiently large $\lambda$, $F$ is topologically conjugate to the full shift on
three symbols, as pictured in Figure~\ref{fig:cubicheur}.
\end{enumerate}

 This proof is similar to the proof for large-scale perturbations of the
 quadratic map so we emphasize the differences.

{\bf Proof of 1.} We prove this in:

(a) Note that \begin{eqnarray*}
\frac{\partial F}{\partial x} (\lambda_0,x)&=&3 x^2 - \lambda_0 + \frac{\partial g}{\partial x} (\lambda_0,x)\\
 &>& 3 x^2 - \lambda_0 - \beta |x| \\
&\ge& -\beta^2/12-\lambda_0.
\end{eqnarray*}
(The last inequality comes from minimizing the quadratic.) Thus as long as $\lambda_0< -\beta^2/12-1$, $F(\lambda_0,x)$ is increasing for all $x$. Since $F(\lambda,x) \to \infty$ for $x \to \infty$ and $F(\lambda,x) \to -\infty$ for $x \to -\infty$, there is a unique fixed point.

\bigskip

{\bf Proof of 2.} Let $s=\sqrt{\lambda_1}$
and $J=[-2s,2s]$. Let $A_1=[-2s,-2.5s/3], A_2=[-s/3,s/3],A_3=[2.5s/3,2s]$, and $A=A_1 \cup A_2 \cup A_3$. Note that the roots of $x^3-\lambda x$ are
$Z=\left\{ 0,\pm \sqrt{\lambda} \right\}$.

\bigskip

(b) Bounds on the derivative. Let $\lambda=\lambda_1$. If $x \in A$,
then $|\partial F/\partial x| \ge 1$.

For $A_2$: If $|x| \le s/3$, then
\begin{eqnarray*}
|\frac{\partial F}{\partial x}|
&\ge& |3 x^2 - \lambda_1| - \beta |x| \ge ( s^2 - \frac{s^2}{3} ) - \frac{\beta s}{3}\\
& \ge & \frac{2s^2 - \beta s}{3}.
\end{eqnarray*}
Thus as long as $\lambda_1> (\beta+\sqrt{\beta^2+24})/4$, $|\partial F /\partial x|> 1$.

For $A_1$ and $A_3$: If $|x|>2.5 s/3 $, $|\partial F/\partial x|>1$:
\begin{eqnarray*}
|\partial F/\partial x|
&\ge& |3 x^2 - \lambda_1| - \beta |x|
\end{eqnarray*}
but for $|x|>2.5s/3$, $ |3 x^2 - \lambda_1| - \beta |x|=3x^2-s^2-\beta x$ is increasing. Thus
\begin{eqnarray*}
|\partial F/\partial x|
&\ge& 3(2.5 s/3)^2-s^2 -\beta (2.5 s/3) \\
&\ge& (3.25 s^2-2.5 \beta s)/3\\
&=& (3.25 \lambda_1 - 2.5 \beta \sqrt{\lambda_1})/3.
\end{eqnarray*}
This last expression is increasing without bound for large $\lambda_1$. Thus as long as $3.25 \lambda_1 -2.5 \beta \sqrt{\lambda_1} -3 >0$, $|\partial F/\partial x|>1$.

(c) If $\lambda=\lambda_1$ and $x \in J \setminus A$, then
$F(\lambda_1,x)$ does not contain points of $J$. Note that the conditions on $g$ guarantee that $|g| < \beta (1 + x^2)$.

Assume that $s/3< |x|< 2.5s/3 $. For $C(x)=x^3-s^2 x$, we have \[|C(s/3)|=8s^3/27 > |C(2.5s/3)|=6.875s^3/27, \] which
implies that over the interval, $\min |C(x)| = 6.875s^3/27$. Thus

\begin{eqnarray*}
|F(\lambda_1,x)|-2s &\ge& 6.875s^3/27 - \beta (1+ 6.25 s^2/9)-2s\\
&\ge& .25 s^3 -.69 \beta s^2 - \beta -2s = .25 \lambda_1^{3/2} - .69 \beta \lambda_1 -2 \sqrt{\lambda_1}- \beta.
\end{eqnarray*}
This expression is 
growing without bound as $\lambda_1 \to \infty$. Thus we can pick $\lambda_1$ sufficiently large that $|F(\lambda_1,x)|-2s >0$ in this region.

(d) If $\lambda \le \lambda_1$, all orbits of $F(\lambda,\cdot)$ are contained in an interval $B$ contained in the interior of $J$:

Let $x>2s$. Then

\begin{eqnarray*}
F(\lambda,x)-x &\ge& x^3-\lambda_1 x - \beta (1 + x^2)-x\\
&\ge& x((x^2-\beta x -(1+\lambda_1)) - \beta \\
& \ge & 2s (4s^2 - \beta (2s) -(1+\lambda_1)) - \beta \\
& = & 8 \lambda_1^{3/2} - 4 \beta \lambda_1 - 2 \sqrt{\lambda_1}(1+\lambda_1) - \beta\\
&= & 6 \lambda_1^{3/2} - 4 \beta \lambda_1 -2 \sqrt{\lambda_1} -\beta.
\end{eqnarray*}
We used the fact that $\lambda_1 > \beta/2$. Note that this last quantity grows arbitrarily large as $\lambda_1 \to \infty$.
Thus as long as $6 \lambda_1^{3/2} - 4 \beta \lambda_1 - 2 \sqrt{\lambda_1} -\beta>0$, $F$ is above the diagonal for $x>2s$. A similar calculation shows that $F$ is below the diagonal for $x<-2s$. Thus there are no orbits outside $J$. 

(e) The sign alternates between each of the three regions for which $F(S)$
is in $S$. In fact, each of $A_i$ for $i=1,2,3$ maps diffeomorphically across $J$. The proof uses the same reasoning as in step (e) of the proof of Corollary~\ref{cor:1D-quadratic}.

\bigskip

This completes the proof of 1 and 2. The result follows from these steps as follows:
We have shown that for large $\lambda$, on the invariant set of $J$, a large-scale perturbation of the cubic map is topologically
 conjugate to the full shift on three symbols. Of the three intervals in $J$, the derivative is positive on the left and right, and negative on the central subinterval. Denote the symbols corresponding to the left, middle, and right subintervals respectively by 
$1_L,-1$, and $1_R$, and use the convention that the sign of $1_L$ and $1_R$ is positive. Thus the sign of the derivative is equal to the sign of the corresponding symbol. Thus the sign of the derivative of a period-$k$ orbit corresponding to a length $k$ sequence of symbols is equal to the sign of the product of the $k$ symbols. 
Since $H_3$ is topologically conjugate to the shift on three symbols, with correspondence between the sign of the derivative of $H_3$ and the sign of the derivative of $F$ for large $\lambda$, each nonflip period-$k$ orbit for $H_3$ is in one-to-one correspondence to a nonflip period-$k$ orbit for $F(\lambda_1,\cdot)$, and each of these orbits is contained in a cascade component with stem period $k$.

There is one nonflip fixed point $P$ for $\lambda \to
-\infty$, and it is an exit orbit. Every nonflip orbit for the shift map on three symbols corresponds to an orbit of $F(\lambda_1,\cdot)$. In other words, 
$IN$ consists of infinitely many orbits, and $OUT$ consists of one orbit, a fixed point. This and Steps (a)-(e) are combined exactly as in the proof of Corollary~\ref{cor:1D-quadratic} to show that by 
Theorem~\ref{thm:main} part $C_0$ all but
possibly one of the nonflip orbits for $\lambda=\lambda_1$ are contained in unique unbounded cascades components.

This completes the proof.
\end{proof}

{\bf The number of nonflip orbits for the three-flip for odd prime $p$.}
The number of cascades of odd prime period $p$
 for a large-scale perturbation of the cubic parametrized map is bounded below by
\[ \frac{(3^p-3)}{2p}.\]
To see this, let $k$ be an odd number.
The number of (not necessarily least) period-$k$
 points is $3^k$, where for the moment we let $k$ be an arbitrary positive integer.
Thus an orbit of period $k$ with nonzero
 orbit index is associated with a symbol sequence with an even number
 of $-1$'s. A simple inductive argument, beginning with $k = 1$, shows that of the $3^k$
 sequences of $k$ symbols, there are $(3^k-1)/2$ with an odd number
 of $-1$'s and $(3^k+1)/2$ with an even number of $-1$'s. Hence all but three sequences correspond to a period-$k$ orbit. Thus there are
 $(3^k+1)/2$ periodic points of (not necessarily least) period $k$
 with nonzero orbit index. 

Now restrict $k=p$ an odd prime. Then the
only period that divides p, but is not equal to p is one. 
Thus there are two nonflip fixed points, which we
 subtract, giving us the formula above.

\section{High-dimensional systems}\label{sec:coupled}

A standard way to get a high dimensional system is to couple low dimensional systems with known properties. We describe a system of $N$ quadratic maps plus coupling. The system has an infinite number of cascades.

In the case of a large-scale perturbation of the quadratic parametrized map, we used the two-symbol $\Gamma$ function to describe the set of nonflip orbits for large $\lambda$. In order to describe the nonflip orbits for the $N$-dimensional case, we now define the more general $\Gamma$ function given in the following definition.

\begin{definition}[Number of periodic points $\Gamma$]\label{D:Gamma}
 For positive integer $N$ and $x = (x_1, \cdots, x_N) \in R^N $, let
$T_N: R^N \to R^N $ be the product of $N$ tent maps,

\[T_N(x) = (T(x_1), \cdots, T(x_N)).\]

For each $k$, let $\Gamma(N,k)$ denote the number of nonflip orbits of period $k$ for
 $T_N$.
\end{definition}
The definition of $\Gamma(N,k)$ for $N=1$ in the introduction is a special case of this definition. The values of $\Gamma(N,k)$ are related, but not in a straightforward way to the number of period-$k$ orbits for the shift map on $2^N$ symbols.

\begin{corollary}[Systems of coupled quadratic parametrized maps]\label{cor:QuadSys}
Let $F:R \times R^N \to R^N $ and $g:R \times R^N \to R^N $ be smooth, and
let each component $F_i$
(for each $i \in 1, 2 \dots, N$)
have the form
\[F_i(\lambda,x_1,\dots,x_N)= K_i(\lambda)-x_i^2+
g_i(\lambda,x_1,\dots,x_N),\] where $g$ is such that for some $\beta>0$, $||g(\lambda,0)||<\beta$,
$||D_xg(\lambda,x)||<\beta$, $\lim_{\lambda \to \infty} K_i(\lambda)
= \infty$, and
$\lim_{\lambda \to - \infty} K_i(\lambda)= -\infty$. 
Then for a residual set of $g$, for each positive integer $k$, the number
of unbounded cascades with stem period $k$ is $\Gamma(N,k)$.
\end{corollary}

\begin{proof}

Assume that $g$ is such that $F$ satisfies genericity hypothesis~\ref{hyp:2}.

The proof of the corollary mimics the method of proof of
Corollary~\ref{cor:1D-quadratic},
with the change that for each $i$,
$K_i(\lambda)$ plays the role of $\lambda$.
In addition we only count the number of cascades with a given stem period.

We say a map $G$ on $ R^N $ is {\bf expansive} on S if $||DG^{-1}|| < 1$
for all $x \in S$, where $\| \cdot \|$ denotes the operator norm.

In what follows we consider $i$ such that $1\le i \le N$.
Assume that $g$ is in the residual set such that $F$ satisfies genericity Hypothesis~\ref{hyp:2}. The proof of this corollary follows from 
Theorem~\ref{thm:main} part $C_0$ by 
showing that the following facts hold:

\begin{enumerate}
\item For all sufficiently small $\lambda$, there are no orbits.
\item For sufficiently large $\lambda_1$, the maximal bounded invariant set of $F(\lambda_1, \cdot)$ is
conjugate to a shift on $2^N$ symbols. On this set, $F$ is diagonally dominant, and the map is expansive.
\end{enumerate}

We now show that these facts hold.
For $\lambda_1$ sufficiently large that $K_i(\lambda_1) > 0$, let
\begin{eqnarray*}
s_i = \sqrt{K_i (\lambda_1)},
J_i = [- 2 s_i, 2 s_i], \mbox{ and }
J = J_1 \times \dots \times J_N.
\end{eqnarray*}
For $j = 1, \dots , 2^N$, let $A_j$ denote one of the $2^N$ connected components of $A \subset J$ where $A=$
$\left\{x: \frac{s_i}{2} \le |x_i| \le 2 s_i \mbox{ for all } i \right\}$.

We show the following estimates for $\lambda_0$ chosen sufficiently small and $\lambda_1$ chosen sufficiently large:

\begin{enumerate}
\item[(a)] At $\lambda_0$, there are no orbits, and orbits of $F$ are strictly monotonically decreasing in each coordinate.
\item[(b)] At $\lambda_1$ and for any point in $A$, $D_xF$ is diagonally dominant, and the dynamics is expansive.
\item[(c)] At $\lambda_1$, the set $F(\lambda_1, R^N \setminus A)$ contains no points of $J$.
\item[(d)] There is a ball $B \subset R^N $ such that for all $\lambda<\lambda_1$, each bounded trajectory is in $B$.

\item[(e)] At $\lambda_1$, $F(\lambda_1, A_i)$ maps diffeomorphically across $J$, for all $i$.
\end{enumerate}

The proof of 1 and 2 follows from these steps: from (b)-(e), it follows that at $\lambda_1$, the map $F$ on the maximal invariant set in $J$ is
topologically conjugate to the full shift on $2^N$ symbols associated
with the $2^N$ rectangles $A_i$. Specifically, all orbits are
hyperbolic. It follows from (d) that there are no orbits on or near
the boundary of $J$ for all $\lambda$ between $\lambda_0$ and
$\lambda_1$. It follows from (a) that there are no orbits on or
near $\lambda=\lambda_0$. Let $C$ be an arc in $PO_{nonflip}(F)$
in region 
$U= [\lambda_0,\lambda_1] \times \mbox{int} (J)$ 
with index orientation $\sigma: I \to C$
(where $I$ is the interval $[0,1]$). If either of $A=\sigma(0)$ or $Z=\sigma (1)$ is on
the boundary of $U$, then it must be on the piece of the boundary
$\lambda=\lambda_1$, since there are no orbits near any other portion of the boundary.
All orbits on $\lambda=\lambda_1$ are hyperbolic orbits, implying that all orbits near the boundary have the same orbit index as the limit orbit on the boundary. Furthermore, they are $N$-dimensionally unstable orbits ($dim_u = N$), implying that they are either flip orbits or that the orbit index is equal to
$(-1)^N$. Therefore, if $N$ is even, there are no entry orbits for $U$, and if $N$ is odd, there are no exit orbits for $U$. This implies that when $N$ is even, $IN$ is empty and $OUT$ is nonempty. When $N$ is odd, $IN$ is nonempty and $OUT$ is empty.
Furthermore, the sign of the derivatives of $F(\lambda_1,\cdot)$ and $T_N$ are the same, implying that there is a one-to-one correspondence between the nonflip period-$k$ orbits of $T_N$ and of $F(\lambda_1,\cdot)$.
By
Theorem~\ref{thm:main} part $C_0$
every nonflip orbit of $F(\lambda_1,\cdot)$ is in a unique unbounded cascade component of stem period $k$.

To complete the proof of this corollary, it only remains to verify the five estimates listed above.

\begin{enumerate}

\item[Proof of (a)] This follows from Corollary~\ref{cor:1D-quadratic}. There
 exists a parameter value $\lambda_0$ sufficiently negative that for any
 values of $x_k$ for $k \ne i$, there are no orbits for
 $K_i(\lambda_0)-x_i^2+g_i(\lambda_0,x_1,\dots,x_n)$. Thus the system
 $F$ has no orbits. In particular, each coordinate of each trajectory is
strictly monotonically decreasing. (It is sufficient to establish this for one coordinate.)

\item[Proof of (b)] If both $x$ and $F(\lambda_1, x)$ are in $J$, and
 $\lambda_1$ is much larger than $g(\lambda_1,x)$, then each
 $K_i(\lambda_1)$ is close in size to $x_i^2$. Hence $\partial
 F_i/\partial x_i$ is large, arbitrarily large as we increase
 $\lambda_1$, while the partial derivatives of $g$ are bounded. Hence
 $D_xF$ is diagonally dominant for large $\lambda_1$. In addition,
 since the diagonal entries grow arbitrarily large, $F$ is expansive
 on $A$.

\item[Proof of (c)] By the assumptions on $g$ and the equivalence of norms in finite dimensions, there is a $\beta_1$ such
 that $|g_i(\lambda,0)|<\beta_1$ for every $i$, and $|\frac{\partial
 g_j}{ \partial x_i}(\lambda,x)| < \beta_1$ for every $i$ and $j$.
For any point in $J \setminus A$, there is an $i$ such that
 $|x_i| < s_i/2$. Therefore, by
Corollary~\ref{cor:1D-quadratic}, for sufficiently large $\lambda_1$, if $|x_i| \le
 s_i/2$, $F_i(\lambda_i,x_i) - 2s_i > 0$, implying that the image of $x$ is outside of $J$.
If $x \in R^N \setminus J$, then for some $i$, $|x_i|>2s_i$. 
By Corollary~\ref{cor:1D-quadratic},
for sufficiently large $\lambda_1$, $F_i(\lambda_1,x)+2s_i < 0$. Thus $F(\lambda_1,x)$ is not in $J$.

\item[Proof of (d)] Let \[ M=\sup_{i \in \left\{1, \dots, N \right\} ,\lambda<\lambda_1}K_i(\lambda).\] $M$ is
 finite, since each $K_i$ is continuous and limits to $-\infty$ as
 $\lambda \to -\infty$. We assume that $M>2$. (If not, set $M=2$.) Let $B=\left\{ x: |x_i|<M+ 2 \beta_1 +1 \mbox{ for } i = 1, \dots, N
 \right\}$, where $\beta_1$ is defined in the proof of (c). Then for all $i$, $|g_i(\lambda,x)|<\beta_1 (1+|x|)$.
For every $x \in R^N \setminus B$, there exists $i$ such
 that $|x_i|>M+ 2 \beta_1+1$. Thus $F_i(\lambda,x)-x_i< M+\beta_1+(\beta_1+1)|x_i| -
 x_i^2<(M+\beta_1)(1-|x_i|)<-(M+\beta_1)$. Therefore the $i^{th}$ component of the orbit is strictly monotonically 
 decreasing by a positive amount, implying that the trajectory is not bounded.

\item[Proof of (e)] Let $A_j$ be a component of $A$. We first show that if $\lambda_1$ is such that for all $i$,
 $\beta<2 s_i$, then there is a point $z \in A_j$ which maps into
 $J$. Specifically, choose $z \in A_j$ such that for all $i$,
 $|z_i|=s_i$. Then for each $i$, $|F_i(\lambda_1,z)| \le
 |g_i(\lambda_1,z)| \le \beta < 2 s_i$.

Fix $A_j$. Let $\lambda_1$ be large enough for (b)-(d) to hold. Let $C = F(\lambda_1, A_j)$. We now show that $J
\subset C$, and $F$ maps $A_j$ diffeomorphically on $J$. In (c) and (d),
we have shown that the image of the boundary of $A_j$ does not
intersect $J$. Since $D_xF$ is nonsingular on $A_j$, $F$ is a local
diffeomorphism on $A_j$. Therefore the interior of $A_j$ maps to an
open set. The set $C$ is compact, so each point of the boundary of $C$
is the image of a boundary point of $A_j$.

As shown above, there exists a point $z \in J \cap F(A_j)$.
Let $y$ be any point not in the compact set $C$. On the straight line
segment from $y$ to $z$, let $w$ be the point closest to $y$ that is in
$C$. Then $w$ is in the image of the boundary of $A_j$, so $w$ is not in
$J$. Since $z$ is in the convex set $J$ but $w$ is not, $y$ is not in $J$.
Hence no point outside of $C$ is in $J$. That is, $J$ is a subset
of $C$.

\end{enumerate}

This completes the proof of the five estimates, thus completing the proof of the corollary.
\end{proof}

\section{Acknowledgements}

Thank you to Safa Motesharrei for his corrections and comments. E.S. was partially supported by NSF Grant DMS-0639300
and NIH Grant R01-MH79502. J.A.Y. was partially supported by NSF Grant DMS-0616585 and NIH Grant
R01-HG0294501.

\begin{flushleft}
%
% References
%
\addcontentsline{toc}{subsection}{References}
\footnotesize
%
% Add in the bbl-file the command \parskip 0pt.
%

{\bf AMS Subject Classification: 37.}\\[2ex]

{\bf Keywords:} Bifurcation, cascades, period doubling, orbit index, horseshoe. \\[2ex]

% Write more than one author separately if they have different
% affiliations, otherwise write the names on the same line, separated
% by commas.
E.~Sander
Department of Mathematical Sciences,
George Mason University,
4400 University Dr.,
Fairfax, VA, 22030, USA.
E-mail: \texttt{esander@gmu.edu}

J.A.~Yorke
Department of Mathematics, IPST, and Physics Department,
University of Maryland,
College Park, MD 20742, USA.
E-mail: \texttt{yorke@umd.edu}

\end{flushleft}

\end{document}